\documentclass[a4paper,10pt]{amsart}

\usepackage[centertags]{amsmath}
\usepackage{amsfonts}
\usepackage{amssymb}
\usepackage{amsthm}
\usepackage{newlfont}
\usepackage{euscript}

\newcommand{\1}{1\!\!{\mathrm I}}
\renewcommand{\Re}{{\Bbb R}}
\newcommand{\eps}{\varepsilon}
\newcommand{\kap}{\varkappa}

\newcommand{\bx}{{x_*}}

\newcommand{\bt}{t_*}

\newcommand{\ax}{\Re^+}

\newcommand{\Es}{\mathsf{E}}
\newcommand{\Cs}{\mathsf{C}}
\newcommand{\Ps}{\mathsf{P}}
\newcommand{\Xs}{\mathsf{X}}

\newcommand{\Ff}{{\EuScript F}}

\newcommand{\Ef}{{\EuScript E}}
\newcommand{\Bf}{{\EuScript B}}
\newcommand{\Tf}{{\EuScript T}}
\newcommand{\Cf}{{\EuScript C}}

\newcommand{\Df}{{\EuScript D}}

\newcommand{\Gf}{{\EuScript G}}
\newcommand{\Qf}{{\EuScript Q}}
\newcommand{\Af}{{\EuScript A}}

\newcommand{\Jf}{{\EuScript J}}

\newcommand{\Pf}{{\EuScript P}}
\newcommand{\pf}{\Pi_{fin}}

\newcommand{\ZZ}{{\Bbb Z}}
\newcommand{\NN}{{\Bbb N}}
\newcommand{\QQ}{{\Bbb Q}}

\newcommand{\DD}{{\Bbb D}}

\newcommand{\demo}{\emph{Proof.} }

\newcommand{\Span}{\mathrm{span}\,}
\newcommand{\supp}{\mathrm{supp}\,}

\newcommand{\bal}{\hbox{{\boldmath $\alpha$}}}
\newcommand{\bbe}{\hbox{{\boldmath $\gamma$}}}

\newcommand{\bfi}{\hbox{{\boldmath $\phi$}}}
\newcommand{\bpsi}{\hbox{{\boldmath $\psi$}}}
\newcommand{\beps}{\eps_*}

\newcommand{\be}{\begin{equation}}
\newcommand{\ee}{\end{equation}}
\newcommand{\nvar}[1]{\left\|#1\right\|_{var}}
\newcommand{\bmu}{{\mu_{inv}}}
\newcommand{\tov}{\mathop{\longrightarrow}\limits^{var}}
\newcommand{\eqd}{\mathop{=}\limits^{df}}
\newcommand{\tN}{$\mathbf{N}$ }

\newcommand{\tR}{$\mathbf{R}$ }
\newcommand{\tLD}{$\mathbf{LD}$ }
\newcommand{\tA}{$\mathbf{A}$}
\newcommand{\tB}{$\mathbf{B}$}
\begin{document}

\numberwithin{equation}{section}
\theoremstyle{plain}
\newtheorem{thm}{Theorem}[section]
\newtheorem*{thm*}{Theorem}
\newtheorem{cor}[thm]{Corollary}
\newtheorem{lem}[thm]{Lemma}
\newtheorem{prop}[thm]{Proposition}
\theoremstyle{definition}
\newtheorem{dfn}[thm]{Definition}
\newtheorem{ex}[thm]{Example}
\newtheorem*{rem}{Remark}

\title[Exponential ergodicity of the solutions to SDE's with a jump noise]
{Exponential ergodicity of the solutions to SDE's with a jump
noise}
\author{Alexey M.Kulik}%
\address{Institute of Mathematics,
Ukrai\-ni\-an National Academy of Sciences, 3, Tereshchenkivska Str., Kyiv 01601, Ukraine}
 \abstract{The mild sufficient conditions for exponential ergodicity of a Markov process, defined
 as the solution to SDE with a jump noise, are given. These
 conditions include three principal claims: recurrence condition
 \textbf{R}, topological irreducibility condition \textbf{S} and
 non-degeneracy condition \textbf{N}, the latter formulated in the
 terms of a certain random subspace of $\Re^m$, associated with the initial
 equation. The examples are given, showing that, in
general, none of three principal claims can be removed without
losing ergodicity of the process. The key point in the approach,
developed in the paper, is that the \emph{local Doeblin condition}
can be derived from  \textbf{N} and \textbf{S} via the
stratification method and criterium for the convergence in
variations of the family of induced measures on $\Re^m$.
 }
\endabstract
\email{kulik@imath.kiev.ua}%
 \maketitle

\centerline{\textsc{Introduction}}
 \vskip 10pt
\setcounter{section}{0}

 In this paper, we study ergodic properties of a Markov process $X$ in $\Re^m$, given
 by an SDE
\be\label{01} d X(t)=a(X(t))dt+\int_{\|u\|\leq 1} c(X(t-),u)\tilde
\nu(dt,du)+\int_{\|u\|>1} c(X(t-),u) \nu(dt,du).
 \ee
Here, $\nu$ is a Poisson point measure,  $\tilde \nu$ is
correspondent compensated measure, and coefficients $a,c$ are
supposed to satisfy usual conditions sufficient for the strong
solution of (\ref{01}) to exist and be unique. Our aim is to give
sufficient conditions for  exponential ergodicity of (\ref{01}),
that impose  as weak restrictions on the L\'evy measure of the
noise,  as it is possible.

There exists two well developed methods to treat the ergodicity
problem for the discrete time Markov processes, valued in a
locally compact phase space. The first one is based on the
\emph{coupling technique} (see detailed overview in
\cite{kloklov_veret}), the second one uses the notions of
\emph{T-chain} and \emph{petite sets} (see  \cite{Meyn_tweedie},
\cite{cline_pu} and references therein). These methods can be
naturally extended to continuous time case either by making a
procedure of time discretization (like it was made for the
solutions to SDE's with jumps in the recent paper \cite{masuda}),
or by straightforward use of the coupling technique in a
continuous time settings (see
\cite{veretennikov_1},\cite{veretennikov_2} for such kind of a
technique for diffusion processes). Typically, in the methods
mentioned above, the two principal features should be provided:

-- recurrence outside some large ball;

-- regularity of the transition probability in some bounded
domain.

The first feature can be provided in a quite standard way via an
appropriate version of the Lyapunov criterium (condition
\textbf{R} in Theorem \ref{t11} below). The second one is more
intrinsic, and requires some accuracy in the choice both of the
concrete terms, in which such feature is formulated, and of the
conditions on the process, sufficient for such feature to hold
true. We deal with the form of the regularity feature, that is
usually called the \emph{local Doeblin condition}, and is
formulated below.

\textbf{LD.} For every $R>0$, there exists $T=T(R)>0$ such that
$$
\inf_{\|x\|,\|y\|\leq R}\int_{\Re^m}[P_x^T\wedge P_y^T](dz)>0,
$$
where $P_x^t(\cdot)\equiv P(X(t)\in \cdot|X(0)=x)$, and, for any
two probability measures $\mu,\varkappa,$
$$
[\mu\wedge \varkappa](dz)\eqd \min\left[{d\mu\over
d(\mu+\varkappa)}(z),{d\varkappa\over
d(\mu+\varkappa)}(z)\right](\mu+\varkappa)(dz).
$$

The non-trivial question is what is the proper form of the
conditions on the coefficients $a,c$ of the equation (\ref{01})
and the L\'evy measure of the noise, sufficient for the local
Doeblin condition to hold true. In a diffusion settings, standard
strong ellipticity (or, more general, H\"ormander type)
non-degeneracy conditions on the coefficients provide that the
transition probability of the process possesses smooth density
w.r.t. Lebesgue measure, and thus \textbf{LD} holds true
(\cite{veretennikov_1},\cite{veretennikov_2}). In a jump noise
case, one can proceed analogously and claim  the transition
probability of the solution to (\ref{01}) to possess a locally
bounded density  (exactly this claim was used as a basic
assumption in the recent paper \cite{masuda}). However, in the
latter case such kind of a claim appears to be too restrictive;
let us discuss this question in more details.
 Consider, for
simplicity, equation (\ref{01}) with $c(x,u)=u$, i.e. a following
non-linear analogue of the Ornstein-Uhlenbeck equation:
\be\label{02} d X(t)=a(X(t))dt+dU_t, \ee where
$U_t=\int_0^t\int_{\|u\|\leq 1}c(u)u\tilde
\nu(ds,du)+\int_0^t\int_{\|u\|> 1}c(u) \nu(ds,du)$ is a L\'evy
process. There exist two methods to provide the process defined by
(\ref{02}) to possess a bounded (moreover, belonging to the class
$C^\infty$) transition probability density. The first one was
initially
 proposed  by J.Bismut (see
\cite{bismut},\cite{bict_grav_Jac},\cite{leandre},\cite{kom_takeuchi}),
the second one -- by J.Picard (see
\cite{picard},\cite{ishikawa_Kunita}). Both these methods require,
among others, the following condition on  the L\'evy measure $\Pi$
of the process $U$ to hold true:
 \be\label{03}\exists \rho\in(0,2):\quad
\eps^{-\rho}\int_{\|u\|\leq \eps}\|u\|^2\Pi(du)\to \infty,\quad
\eps\to0+. \ee This limitation is not a formal one. It is known
(see \cite{Me_jump_reg}, Theorem 1.4), that if \be \label{04}
\mathop{\lim\inf}_{\eps\to 0+}\left[\eps^{2}\ln\left(1\over
\eps\right)\right]^{-1}\sup_{\|v\|=1}\int_{\Re^m}[|(u,v)|\wedge
\eps]^2\Pi(du)=0, \ee  then the transition probability density, if
exists, does not belong to any $L_{p,loc}(\Re^m),p>1$, and
therefore is not locally bounded (note that (\ref{04}) implies
that (\ref{03}) fails). One can say that  when the intensity of
the jump noise is "sparse near zero"\phantom{} in a sense of
(\ref{04}), the behavior of the density essentially differs from
the diffusion one, and the density either does not exist or is
essentially irregular.

On the other hand, let us formulate a corollary of the  general
ergodicity result, given in Theorems \ref{t11},\ref{t13} below.

\begin{prop}\label{p01} Let $m=1,$ suppose that $a(\cdot)$ is locally Lipschitz on $\Re$ and
$\lim\sup\limits_{|x|\to+\infty} {a(x)\over x}<0$. Suppose that
the L\'evy measure $\Pi$ of the process $U$ satisfies the
following conditions:

(i) there exists $q>0$: $\int_{|u|>1}|u|^q\Pi(du)<+\infty$;

(ii) $\Pi(\Re\backslash \{0\}) \not=0$.

Then the solution to (\ref{01}) is exponentially ergodic, i.e. its
invariant distribution $\mu_{inv}$ exists and is unique, and, for
some positive constant $\Cs$, \be\label{05} \forall x\in \Re,\quad
\|P_x^t-\mu_{inv}\|_{var}=O(\exp [-\Cs t]),\quad t\to +\infty. \ee
\end{prop}

In this statement, the non-degeneracy condition (ii) on the jump
noise, obviously, is the weakest possible one: if it fails, then
(\ref{02}) is an ODE, and (\ref{05}) fails also. We can conclude,
that the proper  conditions on the jump noise, sufficient  to
provide exponential ergodicity of the process, defined by
(\ref{02}), are much milder than the conditions that should be
imposed in order to provide that this process possesses regular
(locally bounded or even locally $L_p$-integrable) transition
probability density.

Our way to prove the local Doeblin condition for the solution to
(\ref{01}) strongly relies on the finite-dimensional criterium for
the convergence in variations of the family of induced measures on
$\Re^m$. This criterium was obtained in
\cite{Pilipenko_convergence_by_variation} (the case $m=1$ was
treated in \cite{davydov}).  Via the \emph{stratification method}
(for the detailed exposition of this topic see
\cite{Dav_Lif_Smor}, Section 2) this criterium can be extended to
any probability space with a measurable group of admissible
transformations ($\Leftrightarrow$ \emph{admissible family}), that
generates measurable stratification of the probability space (for
a details, see Section 2 below). The key point is that the
criterium for the convergence in variations of the family of
induced measures is local in the following sense: such a
convergence holds true, as soon as  the initial probability
measure is restricted to any  set, where the gradient (w.r.t.
given admissible family) of the limiting functional is
non-degenerate. We impose a condition (condition \textbf{N} in
Theorem \ref{t13} below), that implies existence of an admissible
family, such that the solution to (\ref{01}) possesses a gradient
w.r.t. this family, and this gradient is non-degenerate with
positive probability. Under this condition, the (local) criterium
for convergence in variation provides the \emph{local Doeblin
condition in a small ball} (Lemma \ref{l210} below). Together with
\emph{topological irreducibility} of the process (provided, in our
settings, by condition \textbf{S} of Theorem \ref{t13}), this
gives condition \textbf{LD}.

 In our construction, we use time-stretching transformations of the
L\'evy jump noise. Such a choice is not the only possible; for
instance, one can use groups of transformations,  varying values
of the jumps (see  \cite{Dav_Lif_84}), and obtain another version
of the non-degeneracy condition \textbf{N}. In order to make
exposition compact, we do not give exact formulation and proof of
the correspondent statement, although the general scheme is
totally the same. Let us just outline that the main advantage of
choice of the differential structure, made in the present article,
is that, for the L\'evy jump noise, time-stretching
transformations, unlike transformations of the phase variable, are
admissible without any regularity claim on the L\'evy measure of
the noise.

The structure of the paper is the following. In Section 1, we
formulate the main statements of the paper. In Section 2, we give
necessary background from the stochastic calculus involving
time-stretching transformations of the L\'evy jump noise, and its
applications to the convergence in variation of the distributions
of the solutions to SDE's with such noise. In Section 3, we prove
the main statements of the paper. Sufficient conditions for the
basic conditions \textbf{R},\textbf{N},\textbf{S} from Section 2,
easy to deal with, are given in Section 4. In Section 5, we give
counterexamples showing that, in general, none of the basic
conditions can be removed without losing ergodicity of the
process.

\section{The main results}
 Let us introduce notation.  Everywhere below
 $\nu$ is a Poisson
point measure on $\ax\times \Re^d$, $\Pi$ is its L\'evy measure
and $\tilde \nu(dt,du)\equiv \nu(dt,du)-dt\Pi(du)$ is the
correspondent compensated point measure. We denote by $p(\cdot)$
the point process associated with $\nu$ and by $\Df$ the domain of
$p(\cdot)$. Later on, we will impose such a conditions on the
coefficients $a,c$ of the equation (\ref{01}),
 that this equation, endowed by the initial condition
$X(0)=x\in\Re^m,$ has the unique strong solution $\{X(t)\equiv
X(x,t),t\geq 0\}$, that is a process with c\'adl\'ag trajectories.
We denote by $\mathsf{P}_\mu$ the distribution in $\DD(\ax,\Re^m)$
of the solution $X(\cdot)$ to (\ref{01}) with Law($X(0))=\mu$, by
$\Es_\mu$ the expectation w.r.t. $\Ps_\mu$, and by $P_\mu^t$ the
distribution of $X(t)$ w.r.t. $\Ps_\mu$. In particular, we denote
$\Ps_x\equiv \Ps_{\delta_x}, P_x^t\equiv
P^t_{\delta_x},x\in\Re^m$.

All the functions used below are supposed to be measurable
(jointly measurable w.r.t. $(x,u)$, if necessary).  The gradient
w.r.t. the variable $z$ is denoted by $\nabla_z$.  The unit matrix
in $\Re^m$ is denoted by $I_{\Re^m}$. We use the same notation
$\|\cdot\|$ for the Euclidean norms both in $\Re^m$ and $\Re^d$,
and for an appropriate matrix norms. The open ball in $\Re^m$ with
the center $x$ and radius $R$ is denoted by $B_{\Re^m}(x,R)$.  The
space of probability measures on $\Re^m$ is denoted by $\Pf$,
   the total variation norm is denoted by $\nvar{\cdot}$. The (closed) support of the measure
$\mu\in\Pf$ is denoted by $\supp \mu$. For
   $\mu\in \Pf$ and non-negative measurable function $\phi:\Re^m\to \Re$, we denote
   $$\phi(\mu)\equiv
\int_{\Re^m}\phi(x)\mu(dx)\in\ax\cup\{+\infty\}.$$

The coefficient $a$ is supposed to belong to the class $C^1(\Re^m,
\Re^m)$ and to satisfy the linear growth condition. In our
considerations, we will deal with the following two types of SDE's
with a jump noise.

{\bf A. Moderate non-additive noise.} The SDE of the type
(\ref{01}) with the jump coefficient $c$ dependent on space
variable $x$. We claim the following standard conditions to hold
true: \be\label{10} \|c(x,u)-c(y,u)\|\leq K(1+\|u\|)\|x-y\|,\quad
\|c(x,u)\|\leq \psi_*(x)\|u\|, \quad u\in\Re^d, x,y\in \Re^m\ee
with some constant $K\in\ax$  and some function $\psi_*: \Re^d\to
\ax$ satisfying linear growth condition. We also claim the
following specific moment condition: \be\label{11}
\int_{\Re^d}\sup_{\|x\|\leq R}\left[\|c(x,u)\|+\|\nabla_x
c(x,u)\|_{\Re^{m^2}}\right]\, \Pi(du)<+\infty,\quad R\in\ax \ee
(the gradient $\nabla_x c(x,u)$ is supposed to exist, and to be
continuous w.r.t. $x$). We interpret this condition in a sense
that the jump part of the equation is \emph{moderate}.

{\bf B. Arbitrary additive noise.} The SDE of the type (\ref{01})
with the jump coefficient $c$ that does not depend on space
variable $x$: $c(x,u)=c(u)$ and $\|c(u)\|\leq K\|u\|$. No moment
conditions like (\ref{11}) are imposed  on the jump part. In this
case, (\ref{01}) is a non-linear analogue (\ref{02}) of the
Ornstein-Uhlenbeck equation.
 Making a change $\nu(\cdot)\mapsto \nu_c(\cdot),
\nu_c([0,t]\times A)\equiv \nu([0,t]\times c^{-1}(A))$, one can
reduce (\ref{02}) to the same equation with $m=d$ and $c(u)=u$. In
order to simplify notation, in a sequel we consider such an
equations only.

 In the both cases given above,  equation
(\ref{01}), endowed by the initial condition $X(0)=x$, has the
unique strong solution, that is a Feller  Markov process with
c\'adl\'ag trajectories,  and $\{P_x^t(\cdot),t\in\ax,x\in\Re^m\}$
is its transition probability. In the case {\bf A}, the
trajectories of this solution a.s. have bounded variation on every
finite interval.

 Denote, like in  \cite{masuda},
$$\Qf = \{f\in C^2(\Re^m, \Re)|\exists  \tilde f \hbox{ locally bounded such that }
\int_{\|u\|>1} f (x + c(x,u))\Pi(du)\leq \tilde f (x), \quad x\in
\Re^m\},
$$
and, for $f\in \Qf$, write
$$
\Af f(x)= \int_{\Re^d}\left[f(x + c(x, u)) - f (x) - (\nabla f(x),
c(x,u))_{\Re^m}\cdot \1_{\|u\|\leq 1}\right]\Pi(du), \quad
x\in\Re^m.
$$

Let us formulate two general statements concerning   convergence
rate of $P_\mu^t$ to the ergodic distribution of $X$ and estimates
for the $\beta$-mixing coefficients of $X$.

\begin{thm}\label{t11}
Suppose  that the \tLD  condition holds true together with the
following recurrence condition:

$\mathbf{R.}$ There exist  function $\bfi\in \Qf$ and constants
$\bal,\bbe>0$ such that
$$
\Af\bfi\leq -\bal \bfi+\bbe\quad \hbox{ and } \bfi(x)\to +\infty,
\quad \|x\|\to +\infty.
$$

Then the process $X$ possesses unique invariant distribution
$\bmu\in\Pf,$ and there exist  positive constants $\Cs_1,\Cs_2$
such that, for every $\mu\in\Pf$ with $\bfi(\mu)<+\infty$,
\be\label{111} \nvar{P_\mu^t-\bmu}\leq
\Cs_1[\bfi(\mu)+1]\exp\left[-\Cs_2t\right],\quad t\in\ax. \ee
\end{thm}

Recall that the $\beta$-mixing coefficient for $X$ is defined by
$$
\beta_{\mu}(t)\equiv
\sup_{s\in\ax}\Es_\mu\left\|\Ps_\mu(\cdot|\Ff_0^s)-\Ps_\mu(\cdot)
\right\|_{var,\Ff_{t+s}^\infty},\quad t\in\ax,
$$
where $\Ff_a^b\equiv \sigma (X(s),s\in [a,b]),
\Ps_\mu(\cdot|\Ff_0^s)$ denotes the conditional distribution of
$\Ps_\mu$ w.r.t. $\Ff_0^s$, and
$$
\left\|\varkappa\right\|_{var,\Gf}\equiv \sup_{\begin{array}{c}
B_1\cap B_2=\emptyset,B_1\cup B_2=\Re^m\\
B_1,B_2\in\Gf
\end{array}}[\varkappa(B_1)-\varkappa(B_2)].
$$
If $\mu=\bmu$, then $\beta(\cdot)\equiv\beta_\bmu(\cdot)$ is the
mixing coefficient of the  stationary version of the process $X$.

\begin{thm}\label{t12} Suppose conditions of  Theorem \ref{t11}
to hold true. Then
\begin{itemize}
\item[(i)] for every $\mu\in\Pf$ with $\bfi(\mu)<+\infty$
\be\label{12} \beta_\mu(t)\leq \Cs_1[\bfi(\mu)+1]\exp[-\Cs_2t];
\ee

\item[(ii)] $\bfi(\bmu)<+\infty$, and thus the mixing coefficient
$\beta(\cdot)$ allows the exponential estimate (\ref{12}).
\end{itemize}
\end{thm}

In order to shorten exposition, we do not formulate here typical
applications of the estimates of the type (\ref{12}), such as
Central Limit Theorem, referring the reader to the literature
(see, for instance, Theorem 4 \cite{veretennikov_2}).

 The following theorem, that   gives sufficient conditions
for condition \textbf{LD} to hold true, is the main result of the
present paper. We need some additional notation.

In the case \textbf{A}, put $\tilde
a(\cdot)=a(\cdot)-\int_{\|u\|\leq 1} c(\cdot,u)\Pi(du)$, and
 denote
$$\Delta (x,u)=[\tilde a(x+c(x,u))-\tilde a(x)]-[\nabla_x
c(x,u)]\tilde a(x).$$
  In the case \textbf{B}, denote  $\Delta
(x,u)=[a(x+u))-a(x)]$ (note that if in the case \textbf{B}
condition (\ref{11}) holds true, then these two
 formulas define the same function). Denote, by
 $\{\Ef_{s}^t, 0\leq s\leq t\}$,  the solution to the linear SDE
in $\Re^{m\times m}$
$$
\Ef_s^t=I_{\Re^m}+\int_s^t \nabla
a(X(r))\Ef_s^r\,dr+\int_s^t\int_{\|u\|\leq 1}\nabla_x c(X(r-))
\Ef_s^{r-}\tilde \nu(dr,du)+\int_s^t\int_{\|u\|>1}\nabla_x
c(X(r-)) \Ef_s^{r-}\nu(dr,du)
$$
($\Ef_s^t$ is supposed to be right continuous w.r.t. both time
variables $s,t$), and define the random linear space $S_t$ as the
span (in $\Re^m$) of the set
$$
\{\Ef_\tau^t\Delta(X(\tau-), p(\tau)), \tau\in \Df\cap (0,t)\}.
$$

\begin{thm}\label{t13} Let the two following conditions to hold
true.

$\mathbf{N.}$ There exist $\bx\in \Re^m, \bt>0$ such that
$$
\Ps_{\bx}(S_{\bt}=\Re^m)>0.
$$

$\mathbf{S.}$ For any $R>0$ there exists $t=t(R)$ such that
$$
\bx\in \supp P_{x}^{t}, \quad \|x\|\leq R.
$$

Then condition \tLD holds true.
\end{thm}

In Section 4 below, we give some sufficient conditions for
\textbf{R},\textbf{N},\textbf{S} to hold true, formulated in the
terms of the coefficients of the equation (\ref{01}) and L\'evy
measure of the jump noise.

\section{Time-stretching transformations and
convergence in variation of induced measures} In this section, we
give necessary background for the technique involving
time-stretching transformations of the L\'evy point measure, with
its applications to the problem of convergence in variation of the
distributions of the solutions to SDE's with jumps. This technique
is  our tool in the proof of Theorem \ref{t13}. Some of the
statements we give without proofs, referring to the
 recent papers \cite{Me_jump_reg},\cite{Me_conv_var}.

Denote $H=L_2(\Re^+), H_0=L_\infty(\Re^+)\cap L_2(\Re^+),
Jh(\cdot)=\int_0^\cdot h(s)\,ds,h\in H.$ For a fixed $h\in H_0$,
we define  the family $\{T_h^t,t\in\Re\}$ of transformations of
the axis $\Re^+$ by putting $T^t_hx, x\in\ax$ equal to the value
at the point $s=t$ of the solution to the Cauchy problem
\be\label{21} z'_{x,h}(s)=Jh(z_{x,h}(s)),\quad s\in \Re, \qquad
z_{x,h}(0)=x.\ee

Denote $T_h\equiv T_h^1$, then $T_{sh}\circ T_{th}=T_{(s+t)h}$
(\cite{Me_jump_reg}). This means that $\Tf_h\equiv \{T_{th}, t\in
\Re\}$ is a one-dimensional group of transformations of the time
axis $\ax$. It follows from the construction that ${d\over dt}
|_{t=0}T_{th}x=Jh(x), x\in\ax.$ We call $T_h$ the \emph{time
stretching transformation} because, for $h\in C(\ax)\cap H_0$, it
can be informally described in the following way: every
infinitesimal segment $dx$ of the time axis should be stretched by
$e^{h(x)}$ times, and then all the  stretched segments should be
glued together, preserving initial order of the segments
(\cite{Me_jump_reg}).

Denote $\Pi_{fin}=
 \{\Gamma\in \Bf(\Re^d),\Pi(\Gamma)<+\infty\}$ and
define, for $h\in H_0, \Gamma\in \pf$, a transformation
$T_h^\Gamma$ of the random measure $\nu$ by
$$
[T_h^\Gamma \nu]([0,t]\times \Delta)=\nu([0,T_{h}t]\times
(\Delta\cap\Gamma))+ \nu([0,t]\times (\Delta\backslash\Gamma))
,\quad t\in\Re^+,\Delta\in \pf.
$$

Further we use the standard terminology from the theory of Poisson
point measures without any additional discussion. The term
"(locally finite) configuration" for a realization of the point
measure is frequently used. We suppose that the basic probability
space $(\Omega,{\Ff},P)$ satisfies condition ${\Ff}=\sigma(\nu)$,
i.e. every random variable is a functional of $\nu$. This means
that in fact one can treat $\Omega$ as the configuration space
over $\ax\times(\Re^d\backslash\{0\})$ with a respective
$\sigma$-algebra.  The image of a configuration of the point
measure $\nu$ under
 $T_h^\Gamma$
 can be described in a following way: every  point
$(\tau,x)$ with $x\not\in \Gamma$ remains unchanged; for every
point $(\tau,x)$ with $x\in \Gamma$, its ``moment of the jump"
$\tau$ is transformed to $T_{-h}\tau$; neither any point of the
configuration is eliminated nor any new point is added to the
configuration.

For $h\in H_0, \Gamma\in\pf$  transformation $T_h^\Gamma$ is
admissible for $\nu$, i.e. the distributions of the point measures
$\nu$ and $T_h^\Gamma \nu$ are equivalent (\cite{Me_jump_reg}).
This imply that (recall that ${\Ff}=\sigma(\nu)$) the
transformation $T_h^\Gamma$ generates the corresponding
transformation of the random variables, we denote it also by
$T_h^\Gamma$.

Define  $\Cf$ as  the set of functionals $f\in
\cap_pL_p(\Omega,P)$ satisfying the following condition: for every
$\Gamma\in \pf$, there exists the random element $\nabla^\Gamma_H
f\in\cap_p L_p(\Omega,P,H)$ such that, for every $h\in H_0$, \be
(\nabla_H^\Gamma f, h)_H=\lim_{\eps\to 0} {1\over \eps}[T_{\eps
h}^\Gamma\circ f-f] \ee with convergence in every $L_p,
p<+\infty$.

Denote
$$ (\rho^\Gamma,h)= -\int_0^\infty
 h(t)\,\tilde\nu(dt,\Gamma), \quad h\in H_0, \Gamma\in \pf.
$$

\begin{lem}\label{l21}(\cite{Me_jump_reg}, Lemma 3.2). For every $\Gamma\in \pf$, the pair
$(\nabla_H^\Gamma,\Cf)$ satisfies  the following conditions:

 1) For every $f_1,\ldots,f_n\in\Cf$ and $F\in C^1_b(\Re^n)$,
$$
F(f_1,\ldots,f_n)\in\Cf\quad \hbox{ and }\quad \nabla_H
F(f_1,\ldots,f_n)=\sum_{k=1}^nF'_k(f_1,\ldots,f_n)\nabla_H f_k
$$
(chain rule).

 2) The map $\rho^\Gamma:h\mapsto (\rho^\Gamma,h)$ is a weak random element
 in $H$ with weak moments of all orders, and
 $$
 E(\nabla_H^\Gamma f,h)_H=-Ef(\rho^\Gamma,h),\quad h\in H,f\in\Cf
$$
(integration-by-parts formula).

3) There exists a countable set $\Cf_0\subset \Cf$ such that
$\sigma(\Cf_0)=\Ff$.
\end{lem}

The construction described before gives us the family
$\Tf=\{T_h^\Gamma, h\in H_0, \Gamma\in\pf\}$ of the {admissible}
transformations of the probability space $(\Omega,\Ff,P)$, such
that the probability $P$ is \emph{logarithmically differentiable}
w.r.t. every $T_h^\Gamma$ with the correspondent logarithmic
derivative equal $(\rho^\Gamma, h)$. This allows us to introduce a
derivative w.r.t. such a family, that is an analogue to the
Malliavin derivative on the Wiener space or the Sobolev derivative
on the finite-dimensional space. However, the structure of the
family $\Tf$ differs from the structure of the family of the
linear shifts: for instance, there exist $h,g\in H_0,\Gamma\in\pf$
such that $T_h^\Gamma\circ T_g^\Gamma\not= T_g^\Gamma\circ
T_h^\Gamma$. This feature motivates the following "refinement",
introduced in \cite{Me_jump_reg}, of the construction described
before.

\begin{dfn}\label{d22}
 A family $\Gf=\{[a_i,b_i)\subset \ax, h_i\in
H_0, \Gamma_{i}\in\pf, i\in\NN\}$ is called  \emph{a differential
grid} (or simply \emph{a grid}) if

(i) for every $i\not= j$,  $\Bigl([a_i,b_i)\times \Gamma_i\Bigr)
\cap \Bigl( [a_j,b_j)\times \Gamma_j\Bigr) =\emptyset$;

(ii) for every $i\in\NN$, $Jh_i>0$ inside $(a_i,b_i)$ and $Jh_i=0$
outside $(a_i,b_i)$.
\end{dfn}

Any grid $\Gf$ generates a partition of some part of the phase
space $\ax\times (\Re^d\backslash \{0\})$ of the random measure
$\nu$ into the cells $\{\Gf_{i}=[a_i,b_i)\times \Gamma_{i}\}$. We
call the grid $\Gf$ \emph{finite}, if $\Gf_i=\emptyset$ for all
indices $i\in\NN$ except some finite number of them. Although
while studying some other problems (such as smoothness of the
transition probability density for $X$, see \cite{Me_jump_reg}) we
typically use infinite grids, in our current exposition we can
restrict ourselves by a finite grids with the number of non-empty
cells equal to $m$ (recall that $m$ is the dimension of the phase
space for the equation (\ref{01})). Thus, everywhere below we
simplify notation from \cite{Me_jump_reg} and consider the grids
$\Gf$ with index $i$ varying from $1$ to $m$.

Denote $T_s^{i}=T_{sh_{i}}^{\Gamma_{i}}$. For any $i\leq m,
s,\tilde s\in \Re $, the transformations $T_s^{i}$,$T_{\tilde s}^{
i}$ commute because so do the time axis transformations $T_{s
h_i}$,$T_{\tilde s h_i}$. Transformation $T_s^{i}$ does not change
points of configuration outside the cell $\Gf_{i}$ and keeps the
points from this cell in it. Therefore, for every $i,\tilde i\leq
m, s,\tilde s\in \Re $, transformations $T_s^{i}$,$T_{\tilde
s}^{\tilde i}$ commute, which implies the following proposition
(\cite{Me_jump_reg}).

\begin{prop}\label{p23} For a given grid $\Gf$ and $t=(t_1,\dots, t_m)\in \Re^m$,
define the transformation
$$
T^{\Gf}_{t}=T^{1}_{t_{1}}\circ T^{2}_{t_{2}}\circ\dots \circ
T^{m}_{t_{m}}.
$$
Then $\Tf^{\Gf}=\{T^\Gf_t, l\in\Re^m\}$ is the group of admissible
transformations of $\Omega$ which is additive in the sense that
$T^\Gf_{t^1+t^2}=T^\Gf_{t^1}\circ T^\Gf_{t^2}, t^{1,2}\in \Re^m.$
\end{prop}

It can be said that, by fixing the grid $\Gf$, we choose from the
whole family of admissible transformations $\{T_h^\Gamma, h\in
H_0,\Gamma\in \pf\}$ the additive sub-family, that is more
convenient to deal with. The following lemma describes the
differential properties of the solution to (\ref{01}) w.r.t. this
family. We denote by $X(x,\cdot)$ the strong solution to
(\ref{01}) with $X(0)=x$.

\begin{lem}\label{l24} $\mathbf{I.}$ In the case \tB,
for every $t\in\ax,\Gamma\in \pf$, every component
$X_k(t),k=1,\dots, m$ of the vector $X(t)$ belongs to the class
$\Cf$. For every $h\in H_0, x\in\Re^m$, the process
$$
Y^{h,\Gamma}(x,t)\equiv ((\nabla_H^\Gamma
X_1(x,t),h)_H,\dots,(\nabla_H^\Gamma X_m(x,t),h)_H)^\top, \quad
x\in \Re^m, t\in\ax
$$
satisfies the equation \be\label{32}
 Y^{h,\Gamma}(x,t)=\int_0^t\int_\Gamma \Delta(X(x,s-),u
 )Jh(s)\,\nu(ds,du)+
\int_0^t [\nabla a](X(x,s)) Y^{h,\Gamma}(x,s)\,ds ,\quad t\geq 0.
\ee

$\mathbf{II.}$ In the case \tA, for every $x\in\Re^m,
t\in\ax,\Gamma\in \pf, h\in H_0$, every component  of the vector
$X(x,t)$ is a.s. differentiable w.r.t.
$\{T_{rh}^\Gamma,r\in\Re\}$, i.e., there exist a.s. limits
\be\label{320} Y_k^{h,\Gamma}(x,t)=\lim_{\eps\to 0} {1\over \eps}
[T_{\eps h}^\Gamma X_k(x,t)-X_k(x,t)],\quad k=1,\dots,m. \ee The
process
$Y^{h,\Gamma}(x,t)=(Y_1^{h,\Gamma}(x,t),\dots,Y_m^{h,\Gamma}(x,t))^\top$
satisfies the equation
$$
 Y^{h,\Gamma}(x,t)=\iint_{[0,t]\times \Gamma} \Delta(X(x,s-),u)Jh(s)\,\nu(ds,du)+
 $$
 \be\label{321}
 +\int_0^t [\nabla a](X(x,s)) Y^{h,\Gamma}(x,s)\,ds +
 \iint_{[0,t]\times \Re^d} [\nabla_x c](X(x,s-),u) Y^{h,\Gamma}(x,s-)\tilde \nu
 (ds,du),\quad t\geq 0.\ee
\end{lem}
Statement \textbf{I} is proved in \cite{Me_jump_reg}, Theorem 4.1;
statement \textbf{II} is proved in \cite{Me_TViMc}, Lemma 4.1.

\begin{rem} Solutions to equations (\ref{32}), (\ref{321}) can be given explicitly:
\be\label{215} Y^{h,\Gamma}(t)=\iint_{[0,t]\times
\Gamma}Jh(s)\cdot
\Ef_s^t\Delta(X(x,s-),u)\nu(ds,du)=\sum_{\tau\in\Df,
p(\tau)\in\Gamma}Jh(\tau)\cdot
\Ef_\tau^t\Delta(X(x,\tau-),p(\tau)). \ee
\end{rem}

For a given grid $\Gf=\{[a_i,b_i)\subset \ax, h_i\in H_0,
\Gamma_{i}\in\pf, i\leq m\}$, denote $Y^{\Gf,i}=Y^{h_i,\Gamma_i},
\quad i=1,\dots,m$ and consider the matrix-valued process
$$
Y^\Gf(x,t)\equiv (Y^i_k(x,t))_{i,k=1}^m,\quad t\in \ax.
$$
The following lemma is the key point in our approach. The
statement of the lemma is formulated  for the cases \textbf{A} and
\textbf{B} simultaneously.

\begin{lem}\label{l25} Let $x\in\Re^m, t>0$ be fixed, denote
$\Omega_{x,t}\equiv\{\det Y^\Gf(x,t)\not=0\}.$ Then \be\label{251}
P|_{\Omega_{x,t}}\circ [X(y,t)]^{-1}\tov P|_{\Omega_{x,t}}\circ
[X(x,t)]^{-1}, \quad y\to x. \ee
\end{lem}

The proof is based on the following criterium for convergence in
variation of induced measures on a finite-dimensional space,
obtained  in \cite{Pilipenko_convergence_by_variation}.

\begin{thm}\label{t_abp}  Let $F, F_n:\Re^m\to \Re^m$ be
measurable functions, that have the approximative derivatives
$\nabla F_n,\nabla F$ a.s. w.r.t Lebesgue measure $\lambda ^m$,
and $E\in {\Bf}(\Re^m)$ has finite Lebesgue measure. Suppose that
$F_n\to F$ and  $\nabla F_n\to \nabla F \quad$ in a sense of
convergence in measure $\lambda^m$ and $\det \nabla F\not=0$ a.s.
on $E$. Then the following statements are equivalent.

(i) for every measurable $A\subset E$ $\quad\lambda^m|_A\circ
F_n^{-1}\mathop{\longrightarrow}\limits^{var} \lambda^m|_A\circ
F^{-1}, n\to +\infty$;

(ii) for every measurable  $A\subset E$ and every  $\delta>0$
there exists a compact set $K_\delta\subset A$ such that
$\lambda^m(A\backslash K_\delta)\leq \delta$ and
$\lim_{n\to+\infty}\lambda^m(F_n(K_\delta))=\lambda^m(F(K_\delta))$.
\end{thm}

In the situation, described in the preamble of the Theorem
\ref{t_abp}, both (i) and (ii) can fail (see, for instance,
Example 1.2 \cite{Me_conv_var}). Thus, in order to provide (i)
(that is our goal), we should impose some additional conditions on
the sequence $\{F_n\}$, sufficient for (ii) to hold true.
 The following two sufficient conditions
were proved in \cite{Pilipenko_convergence_by_variation},
Corollaries 2.5 and 2.7, and in \cite{Me_conv_var}, Theorem 3.1,
correspondingly.

\begin{prop}\label{p27}$\mathbf{I.}$ Let $F_,F_n\in W_{p,loc}^1(\Re^m,
\Re^m)$ with $p\geq m$ ($W_{p,loc}^1$ denotes the local Sobolev
space), and  $F_n\to F,n\to\infty$ w.r.t. Sobolev norm
$\|\cdot\|_{W_p^1(\Re^m,\Re^m)}$ on every ball. Then \be\label{27}
\lambda^m|_A\circ F_n^{-1}\mathop{\longrightarrow}\limits^{var}
\lambda^m|_A\circ F^{-1}, n\to +\infty \hbox{ for every measurable
}A\subset \{\det \nabla F\not=0\}. \ee $\mathbf{II.}$ Let in the
situation, described in the preamble of the Theorem \ref{t_abp},
the sequence $\{F_n\}$ be \textbf{uniformly approximatively
Lipschitz}. This, by definition, means that  for every $\delta>0,
R<+\infty$ there exist a compact set $K_{\delta,R}$ and a constant
$L_{\delta,R}<+\infty$ such that
$\lambda^m(B_{\Re^m}(0,R)\backslash K_\delta)<\delta$
 and every function $F_n|_{K_\delta}$ is a
Lipschitz function with the Lipschitz constant $L_{\delta,R}$. Then (\ref{27}) holds true.
\end{prop}

\emph{Proof of the Lemma \ref{l25}.} Take the sequence $y_n\to x,
n\to +\infty,$ and denote $f=X(x,t), f_n=X(y_n,t)$. It is proved
in \cite{Me_jump_reg} (proof of Theorem 3.1), that the group
$\Tf^\Gf\equiv\{T_s^\Gf, s\in\Re^m\}$  generates a measurable
parametrization of $(\Omega,\Ff,P)$, i.e. there exists a
measurable map $\Phi:\Omega\to \Re^m\times \tilde \Omega$ such
that $\tilde \Omega$ is a Borel measurable space and the image of
every orbit of the group $\Tf^\Gf$ under $\Phi$ has the form
$L\times\{\varpi\}$, where $\varpi\in \tilde \Omega$ and $L$ is a
linear subspace of $\Re^m$. The linear subspace $L$ differs from
$\Re^m$ exactly in the case, when the orbit $T^{\Gf}\{\omega\}$ is
built for such an $\omega$, that, for some $i=1,\dots, m$,
$\nu((a_i,b_i)\times \Gamma_i)=0$ (i.e., some  $T^i$ does not
change $\omega$). For every $\omega$ of such a type $\det
Y^{\Gf}(x,t)=0$, and thus we need to investigate the laws of $f,
f_n,$ restricted to $\Omega_{\Gf}\equiv \{\nu((a_i,b_i)\times
\Gamma_i)>0, i=1,\dots, m\}$, only.

The measure $P|_{\Omega_{\Gf}}$ can be decomposed into a regular
family of conditional distributions such that every conditional
distribution is supported by an orbit of the group $\Tf^\Gf$ (see,
for instance, \cite{Partas}). Therefore we can write \be\label{28}
P(A)=\int_{\tilde \Omega}
P_{\varpi}([A]_{\varpi})\pi(\varpi),\quad A\in \Ff\cap
\Omega_{\Gf}, \ee where $[A]_{\varpi}=\{s\in\Re^m|(s,\varpi)\in
A\}$, $\pi$ is the image of $P|_{\Omega_{\Gf}}$ under the natural
projection $\Omega\to \tilde \Omega$ and $P_\cdot(\cdot)$ is a
\emph{probabilistic kernel}, i.e. $P_\cdot(A)$ is a measurable
function for every $A\in \Bf(\Re^m)$ and $P_{\varpi}(\cdot)$ is a
probability measure on $\Re^m$ for every $\varpi\in \tilde
\Omega$. Any functional $g$ on $\Omega_\Gf$ now can be considered
as a functional on $\Re^m\times \tilde \Omega$, $g=\{g(s,\varpi),
s\in \Re^m,\varpi\in \tilde \Omega\}$. Below, we denote
$[g]_\varpi(\cdot) \equiv g(\cdot,\varpi):\Re^m\to \Re^m$,
$\varpi\in \tilde \Omega$. One can write, for $A\subset
\Omega_{\Gf}$, that
$$
(P|_A\circ g^{-1}) (\cdot)=\int_{\tilde
\Omega}[P_{\varpi}|_{A_\varpi}\circ [g]_\varpi^{-1}](\cdot)
\pi(d\varpi).
$$
Therefore, in order to prove the statement of Lemma \ref{l25}, it is sufficient to prove that
\be\label{281}
P_{\varpi}|_{[\Omega_{x,t}]_\varpi}\circ [f_n]_\varpi\tov P_{\varpi}|_{[\Omega_{x,t}]_\varpi}
\circ [f]_\varpi\quad\hbox{for $\pi$-almost all } \varpi\in\tilde P.
\ee

Denote $\rho_i=(\rho^{\Gamma_i},h_i), i=1,\dots, m,$ and let $\rho^\Gf$  be $\Re^m$-valued
function such that
$(\rho^{\Gf},t)=\sum_{i=1}^mt_i\rho_i, t\in \Re^m$. Then one can show that,
 for $\pi$-almost all $\varpi\in\tilde
\Omega$, the measure  $P_\varpi$ possesses the logarithmic
derivative equal to $[\rho^\Gf]_\varpi$ (we do not give the
detailed exposition here, since this fact is quite analogous to
the one for logarithmically differentiable measures on linear
spaces, see \cite{BogSmol}). One can deduce from the explicit
formula for $\rho^\Gf$ that there exists $c>0$ such that
$E\exp[(\rho^\Gf,t)]<+\infty, \|t\|\leq c$, and, therefore, that
for $\pi$-almost all $\varpi\in \tilde\Omega$ \be\label{210}
\int_{\Re^m}\exp[([\rho^\Gf]_\varpi(s),t)]P_\varpi(ds)<+\infty,
\quad \|t\|\leq c. \ee Due to Proposition 4.3.1 \cite{Boga}, for
every $\varpi$ such that (\ref{210}) holds true, the measure
$P_\varpi$ has the form $P_\varpi(dx)=p_\varpi(x)\lambda^m(dx)$,
where the function $p_\varpi$ is continuous and positive.
Therefore, in order to prove (\ref{281}), it is enough to
 prove that, for $\pi$-almost all such $\varpi$ and every $R>0$,
\be\label{282} \lambda^m|_{B_{\Re^m}(0,R)\cap
[\Omega_{x,t}]_\varpi}\circ [f_n]_\varpi\tov
\lambda^m|_{B_{\Re^m}(0,R)\cap[\Omega_{x,t}]_\varpi} \circ
[f]_\varpi\quad\hbox{for $\pi$-almost all } \varpi\in\tilde P. \ee

In order to prove (\ref{282}), let us introduce auxiliary notions
and give their relations with the notions of Sobolev and
approximative derivatives.

\begin{dfn}\label{d28}
For a given function $F:\Re^m\to \Re^m$ and $\sigma$-finite measure $\kap$ on $\Bf(\Re^m)$
 we say that $F$ is direction-wise $\kap$-a.s. differentiable, if
 there exists function $\nabla F:\Re^{m} \to
 \Re^{m\times m}$ such that, for every $t\in \Re^m$, ${1\over \eps}[F(\cdot+t\eps)-F(\cdot)]\to
 (\nabla F(\cdot),t),\eps\to 0$ $\kap$-a.s. We say that $g$ is direction-wise  differentiable in the
 $L_{p,loc}(\kap)$ sense, if $F\in L_{p,loc}(\kap)$ and,
 for every $t\in \Re^m$,  ${1\over \eps}[F(\cdot+t\eps)-F(\cdot)]\to
 (\nabla F(\cdot),t),\eps\to 0$ in  $L_{p,loc}(\kap)$.
 \end{dfn}

\begin{prop}\label{p28} 1. Let $\kap(dx)=p(x)\lambda^m(dx)$ with $p(x)\geq C_R>0, \|x\|\leq R$ for any
$R>0$. Then every function $F$, that is direction-wise $\kap$-a.s. differentiable, is also
direction-wise $\lambda^m$-a.s. differentiable, and every function $F$, that is direction-wise
differentiable in  $L_{p,loc}(\kap)$ sense, is also
direction-wise differentiable in  $L_{p,loc}(\lambda^m)$ sense. The function $\nabla F$
from the definition of $\kap$-differentiability (either in a.s. or $L_{p,loc}$ sense)
$\lambda^m$-a.s. coincides with the one from the definition of $\lambda^m$-differentiability.

2. If  $F$ is direction-wise $L_{p,loc}(\lambda^m)$-differentiable, then $F\in W_{p,loc}(\Re^m,
\Re^m)$ and $\nabla F$ coincides with its Sobolev derivative.

3. If $F$ is direction-wise $\lambda^m$-a.s. differentiable, then $F$ has approximative derivative
at $\lambda^m$-almost all points $x\in \Re^m$, and $\nabla F$ coincides with its approximative
 derivative.
\end{prop}

\demo Statements 1 and 2 immediately follow from the definition.
Statement 3 follow from Theorem 3.1.4 \cite{Federer} and the
trivial  fact, that the usual differentiability at some point
w.r.t. given direction implies approximative differentiability at
the same point w.r.t. this direction.

Now, we can finish the proof of Lemma \ref{l25}.  Denote
$f_{n,\varpi}=[f_n]_\varpi, f_{\varpi}=[f]_\varpi$. By the
construction, $[T_r^\Gf f_n]_{\varpi}(s)= [f_n]_{\varpi}(s+r),
\quad s,r\in\Re^m, \varpi\in\tilde \Omega$. This, together with
(\ref{28}) and Lemma \ref{l24}, provides that there exists $\tilde
\Omega_0\subset \tilde \Omega$ with
$\pi(\tilde\Omega\backslash\tilde \Omega_0)=0$ such that, for
every $\varpi \in\tilde \Omega_0$, (\ref{210}) holds true, and the
functions $f_{n,\varpi}, f_\varpi$ are either direction-wise
$P_\varpi$-differentiable (in the case \textbf{A}), or
direction-wise differentiable in the $L_{p,loc}(P_\varpi)$ sense
(in the case \textbf{B}), and $\nabla f_{n,\varpi}=
[Y^\Gf(y_n,t)]_\varpi, \nabla f_{\varpi}= [Y^\Gf(x,t)]_\varpi$.
 This, in particular, means that
 $[\Omega_{x,t}]_\varpi=\{\det \nabla f_\varpi\not=0\}$  for $\varpi\in \tilde \Omega_0$.

Now, we can apply the standard theorem on $L_p$-continuity of the
solution to an SDE w.r.t. initial condition (see Theorem 4,
Chapter 4.2 \cite{Gikh_Skor}),
 and obtain that $f_n\to f, [Y^\Gf(y_n,t)]\to [Y^\Gf(x,t)], n\to\infty,$
in $L_p$ sense. This implies that, for $\pi$-almost all $\varpi\in \tilde \Omega_0$,
$f_{n,\varpi}\to f_\varpi, \nabla f_{n,\varpi} \to \nabla f_{\varpi}, n\to \infty$ in
$L_{p}(P_\varpi)$ sense, and, therefore, in $L_{p,loc}(\lambda^m)$ sense. In the case \textbf{B},
convergence (\ref{282}) (and thus the statement of the lemma) follows straightforwardly from
the statement \textbf{I} of Proposition \ref{p27}. In the case \textbf{A}, convergence (\ref{282})
follows from the statement \textbf{II} of the same Proposition, and Lemma 3.3 \cite{Me_conv_var},
that provides that, for $\pi$-almost all $\varpi\in \tilde\Omega,$
 the sequence $\{\nabla f_{n,\varpi}\}$ is uniformly approximatively Lipschitz on $\Re^m$. The lemma
 is proved.

\section{Proofs of the main results}

\subsection{Proof of Theorem \ref{t13}}

We prove Theorem \ref{t13} in two steps. First, we use Lemma
\ref{l25} and  show that, under condition \textbf{N}, the local
Doeblin condition holds true inside some small ball.

\begin{lem}\label{l210} Under condition $\mathbf{N}$,  there exists
$\beps>0$, such that

\be\label{2100} \inf_{x,y\in B(\bx,\beps)}\int_{\Re^m}
\left[P_x^{\bt}\wedge P_y^{\bt}\right](dz)>0. \ee
\end{lem}
\demo Suppose that the grid $\Gf$ is such that, in the notation of
Lemma \ref{l25},
 \be\label{213}
P(\Omega_{\bx,\bt})>0, \ee and denote
$P_{\bx,x}(dz)=\Ps_{x}(\Omega_{\bx,\bt},X(x,\bt)\in dz)$. One can
see that $P_{\bx,x}(dz)=p_{\bx,x}(z)P_x^{\bt}(dz)$ with
$p_{\bx,x}\leq 1$, and therefore
$$
\int_{\Re^m} \left[P_{\bx,x}\wedge P_{\bx,y}\right](dz)\leq
\int_{\Re^m} \left[P_x^{\bt}\wedge P_y^{\bt}\right](dz), \quad x,y
\in\Re^m.
$$
 Thus, in order to prove (\ref{2100}), it is enough to prove that
$$\inf_{x,y\in B(\bx,\beps)}\int_{\Re^m}
\left[P_{\bx,x}\wedge P_{\bx,y}\right](dz)>0.
$$
The latter inequality follows from the condition (\ref{213}), Lemma \ref{l25} and relation
$$
\int_{\Re^m}\left[P_{\bx,x}\wedge P_{\bx,y}\right](dz)=
P(\Omega_{\bx,\bt})-{1\over 2}\|P_{\bx, x}(\cdot)-P_{\bx,
y}(\cdot)\|_{var}\to P(\Omega_{\bx,\bt}),\quad x,y\to \bx.
$$
Thus, the only  thing left to show is that, under condition
\textbf{N}, the grid $\Gf$ can be chosen in such a way that
(\ref{213}) holds true. Denote, by $\Jf_m$, the family  of all
rational partitions of $(0,\bt)$ of the length $2m$; any $J\in
\Jf_m$ is the set of the type
$$
J=\{a_1,b_1, a_2, b_2,\dots, a_m,b_m\}, \quad \hbox{with}\quad
0<a_1<b_1<\dots<a_m<b_m<\bt, \quad a_i,b_i\in \QQ,\quad
i=1,\dots,m.
$$
 Denote, for the set $J$ of such
a type and $r\in \NN$,
$$
\Omega_{J,r}=\biggl\{\forall\, i=1,\dots, m \quad
\exists!\,\tau_i\in (a_i,b_i)\cap \Df, \,\|p(\tau_i)\|\geq {1\over
r},\, i=1,\dots,m
$$
$$
\hbox{ and }\Span\{\Ef_{\tau_i}^t\Delta(X(\bx, \tau_i-),
p(\tau_i)), i=1,\dots,m\}=\Re^m\biggr\}.
$$
Then, elementary considerations show that
$$
\Omega_{\bx, \bt}=\bigcup\limits_{J\in \Jf_m, r\in \NN} \Omega_J.
$$
Therefore, under condition \textbf{N}, there exist
$J^*=\{a_1^*,b_1^*,\dots,a_m^*,b_m^*\}$ and $r^*\in \NN$ such that
$P(\Omega_{J^*,r^*})>0$. Let $h_i^*\in H_0, i=1,\dots, m$  be
arbitrary functions such that, for any $i$, $Jh_i^*>0$ inside
$(a_i^*,b_i^*)$
 and $Jh_i^*=0$ outside $(a_i^*,b_i^*)$. Consider the grid $\Gf^*$
with
$$
\Gamma_i=\begin{cases}\{u|\|u\|\geq {1\over r}\},& i\leq m\\
\emptyset,& i>m\end{cases}, \quad a_i=a_i^*,b_i=b_i^*,
h_i=h_i^*,\, i\leq m,\quad  a_i,b_i,h_i \hbox{ are arbitrary for }
i>m.
$$
 Then formula (\ref{215})
 shows that, on the set $\Omega_{J^*,r^*}$,
 $$
 Y^{\Gf^*,i}(\bx,\bt)\equiv
 Y^{h_i,\Gamma_i}(\bx,\bt)=Jh_i(\tau_i)\Ef_{\tau_i}\Delta(X(\bx,\tau_i-),p(\tau_i)),\quad
 i=1,\dots,m.
 $$
 Since $Jh_i(\tau_i)\not=0$ by the construction, the vectors
 $\{Y^{\Gf^*,i}(\bx,\bt)\}$ are linearly independent iff so are the
 vectors $\{\Ef_{\tau_i}\Delta(X(\bx,\tau_i-),p(\tau_i))\}$.
 Thus,
 $$
 P(\Omega_{\bx,\bt})=P(\det Y(\bx,\bt)\not=0)\geq
 P(\Omega_{J^*,r^*})>0,
 $$
 that gives (\ref{213}). The lemma is proved.

The last step in the proof of Theorem \ref{t13}  is to combine
statement of the previous lemma with the
 condition \textbf{S} and show, that the local  Doeblin condition
 holds true in any bounded region of $\Re^m$.

 \begin{lem}\label{l211} Under conditions \tN and $\mathbf{S}$,
 condition \tLD holds true with
 $$
 T(R)=t(R)+\bt,\quad R>0.
 $$
\end{lem}

\demo The process $X$ is a Feller one; this follows, for instance,
from Theorem 4, Chapter 4.2 \cite{Gikh_Skor}. Therefore, the
function $x\mapsto P_x^t(O)$ is lower semicontinuous for any open
set $O$ and any $t>0$. This, together with condition \textbf{S},
provides that, for any $R>0$,
$$
\delta(R)\equiv \inf_{\|x\|\leq R}P_{x}^{t(R)}(B(\bx,\beps))>0.
$$
Denote
$$
\gamma_*=\inf_{x,y\in B(\bx,\beps)}\int_{\Re^m}
\left[P_x^{\bt}\wedge P_y^{\bt}\right](dz)=\inf_{x,y\in
B(\bx,\beps)}\Big[1-{1\over 2}\nvar{P_x^{\bt}-P_y^{\bt}}\Big]>0,
$$
then, for any $x,y\in B(\bx,\beps)$ and any $A\in \Bf(\Re^m)$.
$$
P_x^{\bt}(A)+P_y^{\bt}(\Re^m\backslash A)\leq 2-2\gamma_*.
$$
Take some $R>0$ and denote $T=T(R)=t(R)+\bt$. Take two independent
processes $X^1, X^2$, satisfying equations of the type (\ref{01})
with the independent point measures $\nu^1,\nu^2$ and starting
from the points $x,y$. Then, for $x,y\in B(0,R)$ and any given
$A\in \Bf(\Re^m)$, we can write
$$
P_x^T(A)+P_y^T(\Re^m\backslash A)=E\Big[
P_{X^1(t(R))}^{\bt}(A)+P_{X^2(t(R))}^{\bt}(\Re^m\backslash
A)\Big]\leq
$$
$$
\leq E\Big[2\1_{\{X^1(t(R))\not\in B(\bt,\beps)\}\cup
\{X^2(t(R))\not\in
B(\bt,\beps)\}}+(2-2\gamma_*)\1_{\{X^1(t(R)),\{X^2(t(R))\in
B(\bt,\beps)}]\leq 2-2\gamma_*\delta^2(R).
$$
Therefore,
$$
\inf_{\|x\|,\|y\|\leq R}\int_{\Re^m}[P_x^T\wedge
P_y^T](dz)=1-{1\over 2}\sup_{\|x\|,\|y\|\leq R}\nvar{P_x^T-
P_y^T}\geq \gamma_*\delta^2(R)>0.
$$
This completes the proof of Lemma \ref{l211} and Theorem
\ref{t13}.

\subsection{Proofs of Theorems \ref{t11},\ref{t12}}
Statements, close to those of  of Theorems \ref{t11},\ref{t12},
are well known in different settings, and there exists several
well developed ways to prove such kind of a statements. For
instance, statement of Theorem \ref{t11} can be derived
straightforwardly from Theorems 5.1, 6.1 \cite{Meyn_tweedie_2},
since condition \textbf{LD} provides that, for any
time-discretized process $X^{\Delta}\equiv \{X(k\Delta), k\in
\ZZ_+\}$ (so called \emph{$\Delta$-skeleton chain}), any compact
set is  a petite set. However,  it is difficult to obtain on this
way an explicit expressions (or estimates) for the constants
$\Cs_1,\Cs_2$, involved in the principal estimates
(\ref{111}),(\ref{12}). Therefore, we use another way to prove
(\ref{111}),(\ref{12}), based on the coupling technique. In
general, we follow the scheme of the proof, proposed for diffusion
processes in \cite{veretennikov_1},\cite{veretennikov_2}, but our
construction of the coupling  slightly differs from the one used
there. This allows us to exclude from the construction auxiliary
conditions, such as Harnack inequality used in
\cite{veretennikov_1} or condition (T) used in
\cite{veretennikov_2}, that are unnatural and restrictive in the
context of SDE's with a jump noise.

Let us start with the construction of the coupling used in the
proof. Since various coupling constructions are used widely in the
literature, we restrict our exposition by the sketch of the
construction only, and omit  technical details.
 First, let us give two basic "bricks" of
our construction. Everywhere below we call "coupling" any
$\Re^m\times \Re^m$-valued process $Y=(Y^1,Y^2)$ such that the
laws of $Y^1,Y^2$ coincide with $\Ps_{\mu_1},\Ps_{\mu_2}$ with
some given $\mu_1,\mu_2\in\Pf$.

\emph{1. Simple coupling.} We call  $Y=(Y^1,Y^2)$ a simple
coupling with a starting point $y=(y_1,y_2)\in\Re^{2m}$, if it is
a coupling with $\mu_{1,2}=\delta_{y_{1,2}}$,  and the processes
$Y^1, Y^2$ are

(a) independent, if $y_1\not =y_2$;

(b) equal one to another, if $y_1=y_2$.

In order to show that such process exists one should simply
consider two equations of the type (\ref{01}) with random point
measures $\nu_1,\nu_2$ that are either independent in the case
$y_1\not =y_2$, or equal one to another in the case $y_1=y_2$.

\emph{2. Gluing coupling.} We call $Y=(Y^1,Y^2)$ a gluing coupling
with a starting point $y=(y_1,y_2)\in\Re^{2m}$ and terminal time
$T>0$, if
 it is a coupling with $\mu_{1,2}=\delta_{y_{1,2}}$,  and

(a) $Y^1=Y^2$, if $y_1=y_2$;

(b) $P(Y^1(T)=Y^2(T))=\int_{\Re^m}[P_{y_1}^T\wedge
P_{y_2}^T](dz)$, if $y_1\not =y_2$.

One can show that such process exists in a following way (one need
to consider the case $y_1\not=y_2$ only). First, due to the
standard \emph{Coupling lemma} (also called \emph{Dobrushin
lemma}, see \cite{dobr}), there exists a probability measure
$\kap$ on $\Re^m\times\Re^m$ such  that
$$
\kap(\{z=(z_1,z_2)|z_1=z_2\}) = \int_{\Re^m} [P_{y_1}^T\wedge
P_{y_2}^T](dz).
$$
Then, this measure is considered as a distribution of $Y$ at the
moment $T$, and the distribution of the whole trajectory of $Y$ is
defined by this measure and the family of  conditional
distributions $\{P(Y\in \cdot|Y(T)=z), z\in \Re^m\times \Re^m\}$
(such a construction is correct since $\DD(\ax, \Re^m\times\Re^m)$
is a Borel measurable space). Any conditional distribution $P(Y\in
\cdot|Y(T)=z), z=(z_1,z_2)\in \Re^m\times \Re^m $ can be
constructed, for instance, as the product of the measures
$$
P(X\in\cdot|X(0)=y_1,X(T)=z_1), \quad
P(X\in\cdot|X(0)=y_2,X(T)=z_2).
$$

Both simple and gluing coupling can be constructed simultaneously
on one probability space $(\Omega_*,\Ff_*,P)$ for all
$y=(y_1,y_2)\in\Re^m\times\Re^m$ in a way, that is jointly
measurable in probability variably $\omega$ and space variable
$y$. For the simple coupling this follows from the standard
theorem on a measurable modification (note that, by the
construction, this coupling is continuous in probability w.r.t.
$y$ on the sets $\{y_1\not=y_2\}$ and $\{y_1=y_2\}$). For the
gluing coupling one can verify this using the \emph{lemma about
three random variables} (\cite{Veretennikov}).   Further we denote
both  these couplings with a starting point $y=(y_1,y_2)$ by
$Y^{y_1,y_2}$.

Now, we can describe our construction. We fix $T,R>0$, that will
be defined later. Construct the probability space $(\Omega,\Ff,P)$
as an infinite product
$$
\Omega=\Omega_0\times \mathop{\Xs}\limits_{k=1}^\infty
\Omega_k,\quad \Ff=\Ff_0\otimes \bigotimes_{k=1}^\infty\Ff_k,\quad
P=P_0\times \mathop{\Xs}\limits_{k=1}^\infty P_k,
$$
where, for $k\geq 1$, $(\Omega_k,\Ff_k,P_k)=(\Omega_*,\Ff_*,P_*)$.
Given $\mu_1,\mu_2\in \Pf,$ construct on $(\Omega_0,\Ff_0,P_0)$
two independent $\Re^m$-valued elements $Z^{1,2}$ with
Law$(Z^{1,2})=\mu_{1,2}$. Next, consider the simple coupling
$\{Y^{y_1,y_2}_1, (y_1,y_2)\in\Re^m\times\Re^m\}$ that is defined
on $(\Omega_1,\Ff_1,P_1)$, and consider the process
$Z_1(t)=Y^{y_1,y_2}(t)\Big|_{y_1=Z^1,y^2=Z^2}, t\geq 0$. Denote
$$
Q_1=\inf\{t|\|Z_1^1(t)\|\leq R, \|Z^2_1(t)\|\leq R\}.
$$
Consider the gluing coupling $\{Y^{y_1,y_2}_2,
(y_1,y_2)\in\Re^m\times\Re^m\}$ that is defined on
$(\Omega_2,\Ff_2,P_2)$, and consider the process
$Z_2(t)=Y^{y_1,y_2}(t-Q_1)\Big|_{y_1=Z^1_1(Q_1),y^2=Z^2_1(Q_1)},
t\geq Q_1$. Denote $Q_2=Q_1+T$. Repeat this construction
iteratively: take the next "independent copy" of the simple
coupling, substitute the terminal value $Z_2(Q_2)$ as the starting
point in it, and wait till the random moment $Q_3$ when both its
coordinates appear inside the ball $\{\|x\|\leq R\}$. Then take
the next "independent copy" of the gluing coupling, substitute the
terminal value $Z_3(Q_3)$ as the starting point in it,  wait till
the moment $Q_4=Q_3+T$, and so on. Define the process
$\{Y(t)=(Y^1(t),Y^2(t)),t\in\ax\}$, by  $Y(t)=Z_k(t),
t\in[Q_{k-1}, Q_k)$; below we call this process "switching
coupling". It has the following properties by the construction:

(i) Law$(Y^i(\cdot))=\Ps_{\mu_i}, i=1,2$;

(ii) for any $k\in \NN$,  $Y^1(t)=Y^2(t),t\geq Q_k$ as soon as
$Y^1(Q_k)=Y^2(Q_k)$;

(iii) $P\Big(Y^1(Q_{2k})=Y^2(Q_{2k})\Big|Y^1(Q_{2k-1})\not
=Y^2(Q_{2k-1})\Big)\geq \inf\limits_{\|x\|,\|y\|\leq
R}\int_{\Re^m}[P_{x}^T\wedge P_{y}^T](dz), k\in\NN$.

Denote $k_*=\min\{k|Y^1(Q_{k})=Y^2(Q_{k})\}$ and put
$Q_*=Q_{k_*}$; $Q_*$ is the "gluing moment" for the coordinates of
the switching coupling $Y$.
 Let us give some estimates that,
together with the property (iii) and condition \textbf{LD}, allow
one  to control the tail probabilities for $Q_*$. Everywhere
below, we suppose that $\bfi\geq 0$. This does not restrict
generality, since one can replace  $\bfi$ by $\bfi+C$ with a
properly chosen constant $C$.

 Fix some $c\in(0,1)$ and take $R>0$ such that $\bfi(x)>\max\Big[{\bbe\over
c\bal}, 1\Big]$ for $\|x\|>R$. Denote $L\equiv \inf\{t|
\|X(t)\|\leq R\}$.

\begin{lem}\label{l33} For any $\mu\in\Pf$ with
$\bfi(\mu)<+\infty$,
$$
\bfi(P_\mu^t)\leq \bfi(\mu) \exp[-\bal t]+{\bbe\over
\bal};\leqno(a)
$$
$$
\Ps_\mu(L>t)\leq  \bfi(\mu)\exp[-(1-c)\bal t].\leqno(b)
$$
\end{lem}
\demo Inequality (a) is a standard corollary of the Dynkin formula
and condition $\Af\bfi\leq -\bal \bfi+\bbe$ (see, for instance,
beginning of the proof of Theorem 6.1 \cite{Meyn_tweedie_2}).
Consider, together with the process $X$, the process $\tilde
X(\cdot)=X(\cdot\wedge L)$ (i.e., the process $X$, stopped at the
first moment of its visit to the ball $\{\|x\|\leq R\}$). By the
construction,  its (extended) generator $\tilde \Af$ satisfies,
for $|x|>R$, the condition
$$
\tilde \Af\bfi(x)=\Af\bfi(x)\leq -\bal\bfi(x)+\bbe\leq
-(1-c)\bal\bfi(x).
$$
Then, writing down the relation, analogous to (a), for the process
$\tilde X$, we obtain that
$$
\Es_\mu(\bfi(\tilde X(t))\1_{\|\tilde X(t)\|>R}\leq
\int_{\|x\|>R}\bfi(x)\mu(dx)\cdot \exp[-(1-c)\bal t]\leq
\bfi(\mu)\exp[-(1-c)\bal t].
$$
Since, by the construction, $\bfi(x)\geq 1$ for $\|x\|>R$, this
implies (b). The lemma is proved.

\begin{cor}\label{c34}
 There exists an invariant measure $\mu_*$ for $X$,  such that
$\bfi(\mu_*)\leq {\bbe\over \bal}$.
\end{cor}
\demo Take some $\mu\in \Pf$ with $\bfi(\mu)<+\infty$ and consider
the family of measures $\{\mu^t, t\in\ax\},$
$$
\mu^t\equiv {1\over t}\int_0^tP_\mu^s\,ds
$$
(the so called \emph{Khasminskii's averages}). It follows from (a)
that $\lim\sup_{t}\bfi(\mu^t)\leq {\bbe\over \bal}$, and this,
together with the Fatou lemma, provides that

(i) the family  $\{\mu^t, t\in\ax\}$ possesses some weak partial
limit $\mu_*$ as $t\to +\infty$;

(ii) $\bfi(\mu_*)\leq {\bbe\over \bal}.$

Moreover, the weak partial limit $\mu_*$ is an invariant measure
for $X$ (the proof of this fact is simple and standard, so we omit
the detailed exposition here). This completes the proof.

Estimates of the Lemma \ref{l33} can be extended from the process
$X$ to the coupling $Y=(Y^1,Y^2)$. The only delicate point here is
that $Y$ does not have to be a Markov process, so we need some
accuracy in writing down the analogues of the estimates (a),(b).
Denote
$$
\bpsi(y)\equiv\bfi(y_1)+\bfi(y_2)\hbox{ and }
\|y\|_\infty\equiv\max[\|y_1\|,\|y_2\|], \quad
y=(y_1,y_2)\in\Re^m\times \Re^m,
$$
then $\bpsi(y)\to +\infty,$ $\|y\|_\infty\to +\infty$. Take
$\tilde R>0$ such that $\bpsi(y)>\max\Big[{2\bbe\over c\bal},
1\Big]$ for $\|y\|_\infty>\tilde R$. Denote $\tilde L\equiv
\inf\{t| \|Y(t)\|_\infty \leq \tilde R\}$.

\begin{lem}\label{l35} (a) Let $\mu_1,\mu_2\in\Pf$ with
$\bfi(\mu_1),\bfi(\mu_2)<+\infty$, and $Y$ be an arbitrary
coupling with $\mathrm{Law}\,(Y^{1,2})=\mu_{1,2}$. Then
$$
E\bpsi(Y(t))\leq [\bfi(\mu_1)+\bfi(\mu_2)] \exp[-\bal
t]+{2\bbe\over \bal}.
$$

(b) Let $\mu_1,\mu_2\in\Pf$ with
$\bfi(\mu_1),\bfi(\mu_2)<+\infty$, and $Y$ be the simple coupling
with $\mathrm{Law}\,(Y^{1,2})=\mu_{1,2}$. Then
$$
P(\tilde L>t)\leq  2[\bfi(\mu_1)+\bfi(\mu_2)] \exp[-(1-c)\bal t].
$$
\end{lem}

\demo The first statement follows immediately from Lemma \ref{l33}
and equality
$E\bpsi(Y(t))=\Es_{\mu_1}\bfi(X(t))+\Es_{\mu_2}\bfi(X(t))$. In
order to prove (b), one should consider separately the cases
$Y^1(0)=Y^2(0)$ and $Y^1(0)\not=Y^2(0)$. In the first case, both
coordinates $Y^1,Y^2$ are the same and move like the process $X$,
i.e., the statement (b) of the Lemma \ref{l33} implies that
$$
P(\tilde L>t, Y^1(0)=Y^2(0))\leq [\bfi(\mu_1)+\bfi(\mu_2)]
\exp[-(1-c)\bal t].
$$
In the second case,  the
 the joint dynamics of the coordinates $Y^1,Y^2$ is described  by the
Markov process, and the  generator $\hat \Af$ of this process
satisfies the relation
$$
\hat \Af \bpsi(y)=\Af \bfi(y_1)+\Af\bfi(y_2)\leq -\bal
\bpsi(y)+2\bbe, \quad y=(y_1,y_2)\in\Re^{m\times m}.
$$
Applying the same arguments with those used in the proof of Lemma
\ref{l33}, we obtain that
$$
P(\tilde L>t, Y^1(0)\not =Y^2(0)) \leq [\bfi(\mu_1)+\bfi(\mu_2)]
\exp[-(1-c)\bal t].
$$
This provides the needed estimate. The lemma is proved.


Using the H\"older inequality, we obtain that \be\label{33}
P(Q_*>t)=\sum_{k=0}^\infty P\Big(Q_*>t, t\in(Q_{2k},
Q_{2k+2}]\Big)\leq \sum_{k=0}^\infty \Big[P(Y^1(Q_{2k})\not =
Y^1(Q_{2k}))\Big]^{1\over 2} \Big[P(Q_{2k+2}\geq t)\Big]^{1\over
2}. \ee Denote
$$
\delta(T,R)=\inf_{\|x\|,\|y\|\leq R}\int_{\Re^m}[P_x^T\wedge
P_y^T](dz).
$$
By the construction, the event $\{Y^1(Q_{2k})= Y^1(Q_{2k}\}$ does
not depend on the values of the process $Y$ up to the moment
$Q_{2k-2}$,  and its probability is not less than $(1-\delta
(T,R))$. Thus, we have an estimate \be\label{35} P(Y^1(Q_{2k})\not
= Y^1(Q_{2k}))\leq (1-\delta (T,R))^k. \ee
 Now let, in the construction
of the switching coupling $Y$ above, $R$ to be  taken equal to
$\tilde R$ given prior to Lemma \ref{l35}, and $T=T(\tilde R)$
(see notation in condition \textbf{LD}). Then $\delta (T,R)>0$,
and (\ref{35}) gives an exponential (w.r.t. $k$) estimate for
$P(Y^1(Q_{2k})\not = Y^1(Q_{2k}))$. Next, $P(Q_{2k+2}\geq t)$ can
be estimated in the following way. Denote $\Delta_k=Q_{k+1}
-Q_{k},$ $\Gf^k=\Ff_{Q_{k}}, k\geq 0$, where $\{\Ff_t,t\in \ax\}$
is the filtration generated by $Y$. Then, for every $k\geq 0$,
$\Delta_{2k-1}=T$ and, from the construction of the switching
coupling $Y$ and statement (b) of Lemma \ref{l35}, we have that
\be\label{34} P(\Delta_{2k}>t|\Ff_{2k})\leq
2[\bfi(Y^1(Q_{2k}))+\bfi(Y^2(Q_{2k}))]\exp[-(1-c)\bal t], \quad
k\geq 0. \ee From the statement (a) of Lemma \ref{l35}, applied to
$t=T$ and $\mu_{1,2}=\mathrm{Law}(Y^{1,2}(Q^{2k-1}))$, we obtain
that
$$
 E[\bfi(Y^1(Q_{2k}))+\bfi(Y^2(Q_{2k}))|\Ff_{2k-1}]\leq
{2\bbe\over \bal}+2\sup_{\|x\|\leq R}\bfi(R),\quad k\geq 1,
$$
here we used that $\|T^{1,2}(Q_{2k-1})\|\leq \tilde R.$ This,
together with (\ref{34}), gives that
$$
P(\Delta_{2k}>t|\Ff_{2k-1})\leq \Big[{4\bbe\over
\bal}+4\sup_{\|x\|\leq R}\bfi(R))\Big]\exp[-(1-c)\bal t], \quad
k\geq 1,
$$
and, consequently, $$  E\Big(\exp[{1-c\over
2}\Delta_{2k}]|\Ff_{2k-1}\Big)\leq \Big[{4\bbe\over
\bal}+4\sup_{\|x\|\leq R}\bfi(R))\Big]\int_0^\infty{1-c\over
2}e^{-{1-c\over 2}z}dz=\Big[{4\bbe\over \bal}+4\sup_{\|x\|\leq
R}\bfi(R))\Big], \quad k\geq 1. $$ Analogously, we have that $$
E(\exp[{1-c\over 2}\Delta_{0}])\leq 2[\bfi(\mu_1)+\bfi(\mu_2)].
$$
At last, $\Delta_{2k+1}=T, k\geq 0$. This and two previous
estimates provide that  \be\label{36} E\exp[{1-c\over
2}Q_{2k+2}]=E\prod_{j=0}^{2k+1} \exp[{1-c\over 2}\Delta_j]\leq
2[\bfi(\mu_1)+\bfi(\mu_2)]\cdot \exp[(k+1){1-c\over 2}T]\cdot
\Big[{4\bbe\over \bal}+4\sup_{\|x\|\leq R}\bfi(R))\Big]^k=\ee $$ =
2 \exp[{1-c\over 2}T][\bfi(\mu_1)+\bfi(\mu_2)]\cdot \exp[kD],\quad
k\geq 0,\quad \hbox{ with } \quad D= {1-c\over
2}T+\ln\Big[{4\bbe\over \bal}+4\sup_{\|x\|\leq R}\bfi(R))\Big].
$$
Take $p=\max\Big[1, -{2D\over \ln(1-\delta(T,R))}\Big]$, then, by
H\"older inequality,
$$
 E\exp[{1-c\over
2p}Q_{2k+2}]\leq \max\left\{2 \exp[{1-c\over
2}T][\bfi(\mu_1)+\bfi(\mu_2)], 1\right\}\cdot
(1-\delta(T,R))^{-{k\over 2}},
$$
and, by Chebyshev inequality,
$$
P(Q_{2k+2}\geq t)\leq \exp[-{1-c\over 2p}t]\cdot\max\left\{2
\exp[{1-c\over 2}T][\bfi(\mu_1)+\bfi(\mu_2)], 1\right\}\cdot
(1-\delta(T,R))^{-{k\over 2}}.
$$
This estimate, together with (\ref{33}),(\ref{35}), gives that
$$
P(Q^*>t)\leq \max\left\{2 \exp[{1-c\over
2}T][\bfi(\mu_1)+\bfi(\mu_2)], 1\right\}\cdot
\Big[1-(1-\delta(T,R))^{{1\over 4}}\Big]^{-1}\cdot \exp[-{1-c\over
4p}t], \quad t\in \ax.
$$
Now, we can use standard arguments (see \cite{veretennikov_1}) and
complete the proofs of Theorems \ref{t11}, \ref{t12}. In order to
prove (\ref{111}), consider the switching coupling with
$\mu_1=\mu$ and $\mu_2=\mu_*$ given by Corollary \ref{c34}. Then
$$
\nvar{P_\mu^t-\mu_*}\leq P(Q^*>t)\leq \tilde
\Cs_1[\bfi(\mu)+1]\exp[-\Cs_2 t],\quad t\in \ax
$$
with $\tilde \Cs_1=\max\Big[ {\bbe\over \bal}, 1\Big]\cdot 2
\exp[{1-c\over 2}T]\cdot \Big[1-(1-\delta(T,R))^{{1\over
2}}\Big]^{-1},\Cs_2={1-c\over 4p}.$ Analogously, for a given
$s,t\in \ax$ consider the switching coupling with $\mu_1=\mu$ and
$\mu_2=P^s_\mu$. Then, by statement (a) of Lemma \ref{l33},
$\bfi(\mu_2)\leq \bfi(\mu)+{\bbe\over \bal}$, and
$$
\Es_\mu\left\|\Ps_\mu(\cdot|\Ff_0^s)-\Ps_\mu(\cdot)
\right\|_{var,\Ff_{t+s}^\infty}\leq P(Q^*>t)\leq
\Cs_1[\bfi(\mu)+1]\exp[-\Cs_2 t],\quad t\in \ax
$$
with $\Cs_1=2\tilde \Cs_1$. These two estimates imply
(\ref{111}),(\ref{12}) with $\mu_{inv}=\mu_*$. Theorems
\ref{t11},\ref{t12} are proved.

\section{Sufficient conditions.}

In this section, we give some sufficient conditions for
\textbf{R},\textbf{N},\textbf{S} to hold true.

\subsection{Condition R} There exists a wide range of conditions, that are
sufficient for the recurrence condition \textbf{R}, see
\cite{masuda}, Section 2.3. Here, we give only one condition of
such a type, that is an analogue of Lemma 2.4 \cite{masuda}, but
with a condition (11) of this Lemma replaced by an essentially
weaker one (condition 3 below).

\begin{prop}\label{p41} Suppose that the following conditions hold true.

\noindent  1. There exist $R,\alpha>0$ such that
$$
(a(x),x)_{\Re^m}\leq -\alpha \|x\|^2,\quad \|x\|\geq R.
$$

\noindent 2. There exists $q\in(0,+\infty)$ such that
$$
\int_{\|u\|>1} \|u\|^{q}\Pi(du)<+\infty.
$$

\noindent 3. The function $c$ can be decomposed into a sum
$c=c_1+c_2$ with $c_1, c_2$ such that

 3a. for some  function  $\psi$, ${\psi(x)\over \|x\|}\to 0$ as $\|x\|\to
\infty$,
$$
 \, \|c_1(x,u)\|\leq \psi(x)\|u\|,\quad u\in\Re^d, x\in
\Re^m;
$$
$$
\|x+c_2(x, u)\|\leq \|x\|, \quad x\in \Re^m, \|u\|>1,\quad
c_2(\cdot,u)\equiv 0, \quad \|u\|\leq 1.\leqno
{\phantom{}\,\,\quad3b.}
$$

Then condition \tR holds true.
\end{prop}

\begin{rem} In the case \textbf{B}, condition 3  holds true automatically with $c_1=c, c_2=0$.
\end{rem}

\demo Consider $\bfi\in C^2(\Re^m)$ such that $\bfi(x)=\|x\|^q,
\|x\|\geq R$. Without losing generality we can suppose that the
constant $R$ is chosen in such a way that $\delta_r\equiv
\sup_{\|x\|\geq r}{\psi(x)\over \|x\|}\leq {1\over 2}, r\geq R.$
Then, for $\|x\|\geq 2R$,
$$
\Af\bfi(x)=q(a(x),x)_{\Re^m}\|x\|^{q-2}+\int_{\Re^m}\Bigl[\|x+c(x,u)\|^q-\|x\|^q-
q\1_{\|u\|\leq 1}(c(x,u),x)_{\Re^m}\|x\|^{q-2}\Bigr]\Pi(du)=
$$
$$
=q(a(x),x)_{\Re^m}\|x\|^{q-2}+\int_{\|u\|\leq
1}\Bigl[\|x+c_1(x,u)\|^q-\|x\|^q- q\1_{\|u\|\leq
1}(c_1(x,u),x)_{\Re^m}\|x\|^{q-2}\Bigr]\Pi(du)+
$$
$$
+\int_{\|u\|>1}\Bigl[\|x+c(x,u)\|^q-\|x\|^q\Bigr]\Pi(du).
$$
Due to condition 3b,
$$
\int_{\|u\|>1}\Bigl[\|x+c(x,u)\|^q-\|x\|^q\Bigr]\Pi(du)\leq
\int_{\|u\|>1}\Bigl[\|x+c(x,u)\|^q-\|x+c_2(x,u)\|^q\Bigr]\Pi(du)=
$$
$$
=\int_{\|u\|>1}\Bigl[\|x(u)+c_1(x,u)\|^q-\|x(u)\|^q\Bigr]\Pi(du),
$$
where $x(u)=x+c_2(x,u)$. If $O$ is an open subset of $\Re^m$,
$\{x+sc, s\in[0,1]\}\subset O$ and $\phi\in C^2(O)$, then
\be\label{41}|\phi(x+c)-\phi(x)-(\nabla \phi(x),c)_{\Re^m}|\leq
\|c\|^2\sup_{s\in[0,1]}\|[\nabla^2 \phi](x+sc)\|\ee (the Taylor's
formula). Then, for $r\geq 2R$,
$$
 \int_{\|u\|\leq 1}\Bigl[\|x+c_1(x,u)\|^q-\|x\|^q-
q(c_1(x,u),x)_{\Re^m}\|x\|^{q-2}\Bigr]\Pi(du)\leq
\mathrm{const}\cdot [\|x\|(1-\delta_r)]^{q-2}[\|x\|\delta_r]^{2},
\quad \|x\|\geq r,
$$
here we applied  (\ref{41}) with $\phi(x)=\|x\|^q$,
$O=\{\|x\|>{r(1-\delta_r)}\}$ and $c=c_1(x,u)$.

If $q\geq 1$, then, applying the inequality
$|\phi(x+c)-\phi(x)|\leq \|c\|\sup_{s\in[0,1]}\|[\nabla
\phi](x+sc)\|$, with the same $\phi$, $x=x(u)$ and $c=c_1(x,u)$,
we obtain analogously
$$
\left|\int_{\|u\|>1}\Bigl[\|x(u)+c_1(x,u)\|^q-\|x(u)\|^q\Bigr]\Pi(du)\right|\leq
\mathrm{const}\cdot [\|x\|(1-\delta_r)]^{q-1}[\|x\|\delta_r],
\quad \|x\|\geq r,
$$
here we used that $\|x(u)\|\leq \|x\|$.  If $q<1$, then we apply
inequality $a^q+b^q\geq (a+b)^{q}, a,b\in\ax$ and write
$$
\left|\int_{\|u\|>1}\Bigl[\|x(u)+c_1(x,u)\|^q-\|x(u)\|^q\Bigr]\Pi(du)\right|\leq
\int_{\|u\|>1} \|c_1(x,u)\|^{q}\Pi(du)\leq \mathrm{const}\cdot
[\psi(x)]^q.
$$
Thus,  $$ \Af\bfi(x) \leq -q\alpha \|x\|^q +C
\|x\|^{q}\Bigl\{(1-\delta_r)^{q-2}[\delta_r]^{2}+
(1-\delta_r)^{q-1}\delta_r +
\delta_r^{2q}\Bigr\}=
$$
$$=\bfi(x)\left[-q\alpha+C \Bigl\{(1-\delta_r)^{q-2}[\delta_r]^{2}+
(1-\delta_r)^{q-1}\delta_r + \delta_r^{2q}\Bigr\}\right] ,\quad
\|x\|\geq r$$
 with some constant $C$. Since $\delta_r\to 0, r\to+\infty,$ this gives \textbf{R}.
The proposition is proved.

\subsection{Condition N} In \cite{Me_TViMc}, \cite{Me_jump_reg}
the conditions were given, sufficient for the set $\{S_t=\Re^m\}$
to have probability one. Below we give a more mild version of
these conditions, sufficint for this set to have non-zero
probability. Denote, for $x\in\Re^m$,
$\Theta_x=\{u|I_{\Re^m}+\nabla_x c(x,u)$ is invertible$\}$, and
put
$$
\hat \Delta(x,u)=[I_{\Re^m}+\nabla_x c(x,u)]^{-1}\Delta(x,u),
\quad u\in \Theta_x.
$$
Denote, by $S^m\equiv\{v\in\Re^m|\|v\|_{\Re^m}=1\}$, the unit
sphere in $\Re^m$.

\begin{prop}\label{p42} Suppose that there exists $x_*\in \Re^m$ such that
\be\label{m} \forall \eps>0, v\in S^m\quad \Pi\Big(u\in
\Theta_{x_*}\Big| \,(\hat \Delta(x_*,u), v)_{\Re^m}\not=0,
\|c(x_*,u)\|< \eps\Big)>0. \ee Then condition \tN holds true with
this $x_*$ and arbitrary $t_*>0$.
\end{prop}

\demo We need to prove that, on some $\Omega_0$ with
$\Ps_{x_*}(\Omega_0)>0$, \be\label{42}
\{\Ef_\tau^t\Delta(X(\tau-), p(\tau)), \tau\in \Df\cap
(0,t)\}=\Re^m. \ee Below, the set $\Omega_0$ will be constructed
explicitly, and in particular, on the set $\Omega_0$, the matrix
$\Ef_0^t$ will be non-degenerate (and therefore, for any variables
$\theta, \tau$, $0\leq \theta\leq \tau\leq t$, the matrix
$\Ef_\theta^\tau$ also will be non-degenerate). Then, on this set,
$$
\{\Ef_\tau^t\Delta(X(\tau-), p(\tau)), \tau\in
\Df\}=\Ef_0^t\{\Ef_0^{\tau-}\hat \Delta(X(\tau-), p(\tau)),
\tau\in \Df\cap (0,t)\},
$$
and (\ref{42}) is equivalent to \be\label{43}
\{[\Ef_0^{\tau-}]^{-1}\hat \Delta(X(\tau-), p(\tau)), \tau\in
\Df\cap (0,t)\}=\Re^m. \ee

For $n\geq 1$, consider the set $\Df^n=\{\tau\in \Df,
\|p(\tau)\|\geq {1\over n}\}$. This set is a.s. locally finite,
and therefore can be enumerated increasingly,
$\Df^n=\{\tau_1^n,\tau_2^n,\dots\}$. Denote, for $k\leq m$,  $S_k
^n=\{[\Ef_0^{\tau_j^n-}]^{-1}\hat \Delta(X(\tau_j^n-),
p(\tau_j^n)), j\leq k\}$. By the construction, $S_k^n$ is a linear
span of a finite family  of vectors. Let us consider the $k$-th
vector from this family,
$$
[\Ef_0^{\tau_k^n-}]^{-1}\hat \Delta(X(\tau_k^n-), p(\tau_k^n)).
$$
One can construct the measurable map $V:(\Re^m)^{k-1}\to S^m$ such
that
$$
\forall j=1,\dots, k-1\quad V(x_1,\dots, x_{k-1})\perp x_j, \quad
x_1,\dots, x_{k-1}\in \Re^m.
$$
We will write $v_{k-1}^n\equiv
\Big([\Ef_0^{\tau_j^n-}]^{-1}\Big)^*\cdot
V(\{[\Ef_0^{\tau_j^n-}]^{-1}\hat \Delta(X(\tau_j^n-),
p(\tau_j^n)), j< k\}),$ where $(M)^*$ denotes the adjoint matrix
for $M$. The random vector $v_{k-1}^n$ is well defined on the set
$\{\Ef^{\tau_k^n-}$ is invertible$\}\in \Ff_{\tau_{k^n-}}$, and is
$\Ff_{\tau_{k^n-}}$-measurable.

The value $p(\tau_k^n)$ is independent of $\Ff_{\tau_{k^n-}}$, and
its distribution is equal to ${1\over
\lambda_n}\Pi(\cdot\cap\{\|u\|\geq 1\})$, where
$\lambda_n=\Pi(\|u\|\geq 1)$. Therefore, on the set
$\{\Ef^{\tau_k^n-}$ is invertible$\}$,
$$
P([\Ef_0^{\tau_k^n-}]^{-1}\hat \Delta(X(\tau_k^n-),
p(\tau_k^n))\hbox{ is well defined and }\not \in
S_{k-1}^n|\Ff_{\tau_{k^n-}})=
$$
\be\label{431} ={1\over \lambda_n}\Pi\Big(u\in \Theta_{x}\Big|
\,(\hat \Delta(x,u), v)_{\Re^m}\not=0\Big)\Big|_{x=X(\tau_k^n-),
v=v_{k-1}^n} \geq{1\over \lambda_n}\inf_{v\in S^m}\Pi\Big(u\in
\Theta_{x}\Big| \,(\hat \Delta(x,u),
v)_{\Re^m}\not=0\Big)\Big|_{x=X(\tau_k^n-)}. \ee

For a given $\eps>0$, consider, for $n\geq 1$, the maps
$$
f_n: S^m\ni v\mapsto \Pi\Big(u\in \Theta_{x_*}\Big| \,(\hat
\Delta(x_*,u), v)_{\Re^m}\not=0, \|c(x_*,u)\|< \eps, \|u\|>{1\over
n}\Big).
$$
Since functions $c, \nabla_x c, \hat\Delta$ are continuous w.r.t.
$x$ on their domain, every $f_n$ is lower semicontinuous. For
every $v\in S^m$, $f_n(v)$ monotonously tends to a positive limit
as $n\uparrow \infty$. Therefore, due to Dini theorem, there
exists $n=n(\eps)\in \NN$ such that
$$\inf_{v\in S^m} \Pi\Big(u\in \Theta_{x_*}\Big|
\,(\hat \Delta(x_*,u), v)_{\Re^m}\not=0, \|c(x_*,u)\|< \eps,
\|u\|>{1\over n}\Big)>0.
$$
Analogously, the function
$$
\Re^m \ni x\to \inf_{v\in S^m} \Pi\Big(u\in \Theta_{x}\Big|
\,(\hat \Delta(x,u), v)_{\Re^m}\not=0, \|c(x,u)\|< \eps,
\|u\|>{1\over n(\eps)}\Big)
$$
is lower semicontinuous, and thus there exists
$\delta=\delta(\eps)>0$ such that \be\label{45} \inf_{v\in S^m,
x\in B(x_*, \delta)} \Pi\Big(u\in \Theta_{x}\Big| \,(\hat
\Delta(x,u), v)_{\Re^m}\not=0, \|c(x,u)\|< \eps, \|u\|>{1\over
n(\eps)}\Big)>0. \ee Define iteratively  $\eps_1,\dots, \eps_m$ in
the following way. Take $\eps_m>0$ arbitrary, and put
$$
\eps_{k-1}=\min\Big[{1\over 6}\delta(\eps_k), {1\over 2}
\eps_k\Big],\quad k=2,\dots, m.
$$
By the construction, $\delta(\eps_k)>3\sum_{l<k}\eps_l, k=2,\dots,
m.$ Put $n=\max_{l\leq k} n(\eps_{l})$. Let us show that, for
these  $n\in \NN, $  $\eps_1,\dots, \eps_{m}>0$,  and properly
chosen $r\in(0,{t_*\over m})$, the set
$$
\Omega_0\equiv [\cap_{k=1}^m\Omega_{k}^r]\cap \Omega^r
$$
has non-zero probability, where
$$
\Omega^r=\Big\{\Ef_r^{t_*}\hbox{ is invertible }\Big\},
$$
$$
\Omega_{k}^r=\Big\{\Df^n\cap((k-1) r, kr]=\{\tau_k^r\}, \,
\Ef^{kr}_{(k-1)r}\hbox{ is invertible, }
\|X(\tau_k^n-)-X((k-1)r)\|\leq \eps_{k-1},
\|X(kr)-X(\tau_k^n)\|\leq \eps_k,
$$
$$
\|X(\tau_k^n)-X(\tau_k^n-)\|\leq \eps_k\hbox{ and
}[\Ef_0^{\tau_k^n-}]^{-1}\hat \Delta(X(\tau_k^n-),
p(\tau_k^n)\not\in\mathrm{span}\, \{[\Ef_0^{\tau_j^n-}]^{-1}\hat
\Delta(X(\tau_j^n-), p(\tau_j^n)), j<k\}\Big\}.
$$
It is easy to verify that $P[\Omega^r|\Ff_{mr}]>0$ a.s. (the proof
is omitted), so we need to verify that
$P(\cap_{k=1}^m\Omega_{k}^r)>0.$

The process $X$ can be described in the following way: at the
moments $\tau_k^n,k\geq 1,$ it has the jumps of the value
$c(X(\tau_k^n-), p(\tau_k^n))$, and on every interval of the type
$(\tau_{k-1}^n, \tau_n^k), k\geq 1$, it moves due to SDE
 \be\label{46}
dX^n(t)=a^n(X^n(t))\,dt+\int_{\|u\|< {1\over n}}c(X^n(t-),u)\tilde
\nu(dt,du),  \ee where $a^n(x)=a(x)+\int_{\|u\|\in[{1\over
n},1]}c(x,u)\Pi(du), \tau_0^n=0.$ Denote, by $X^n(x,\cdot)$, the
solution to (\ref{46}) with $X^n(0)=x$ and, by $\Ef^{n,x,
\cdot}_{\cdot}$, the correspondent stochastic exponent (both
$X^n(x,\cdot)$ and $\Ef^{n,x}$ are independent of the point
process $p|_{\Df^n}$). Then, for a given $n\in \NN, \eps_1, \dots,
\eps_m>0$, one can choose $r\in (0,{t_*\over m})$ small enough for
$$
D\equiv \inf_{x\in B(x_*,3\sum_{l=1}^m\eps_m)}P\Big(\forall s\leq
r\, \|X^n(x, s)-x\|<\min_l\eps_l, \Ef_0^{n,x, s}\hbox{ is
invertible}\Big)>0.
$$
Now, let us estimate $P(\Omega^r_k|\cap_{l<k}^m\Omega_{l}^r)$. The
set $\cap_{l<k}^m\Omega_{l}^r$ belongs to $\Ff_{(k-1)r}$, and, on
this set,
$$
\|X((k-1)r)-x_*\|\leq 3\sum_{l<k-1}\eps_l+2\eps_{k-1}.
$$
Now we can take subsequently conditional expectations first w.r.t.
$\Ff_{\tau_k^n-}\vee \sigma(p(\tau_k^n))$, then w.r.t.
$\Ff_{\tau_k^n-}$ (on this step, we use (\ref{431})), and, at
last, w.r.t.
 $\Ff_{(k-1)r}$, and  write that, on
this set,
$$
P\Big(\Omega_k^n\cap\{\|X((\tau_k^n-)-x_*\|\leq
3\sum_{l<k}\eps_l,\|X((\tau_k^n)-X((\tau_k^n-)\|\leq \eps_k,
\|X((kr)-X((\tau_k^n)\|\leq \eps_k \}|\Ff_{(k-1)r}\Big)\geq
$$
$$
\geq D\cdot (r\lambda_n e^{-r\lambda_n})\cdot {\gamma\over
\lambda_n}\cdot D=r\gamma D^2 e^{-r\lambda_n},
$$
where
$$
\gamma=\min_{l\leq k}\inf_{x\in B(x_*, 3\sum_{l<k})}\Pi\Big(u\in
\Theta_{x}\Big| \,(\hat \Delta(x,u), v)_{\Re^m}\not=0, \|c(x,u)\|<
\eps_k, \|u\|>{1\over n}\Big)>0
$$
by the construction. Therefore,
$$
P(\cap_{k=1}^m\Omega_k^r)>0,
$$
that gives the needed statement. The proposition is proved.

Let us also give sufficient condition for \textbf{N}, in which
conditions on the coefficients $a,c$ and L\'evy measure of the
noise are separated.

For $w\in S^d$ (recall that $S^d$ denotes the unit sphere in
$\Re^d$), $\varrho\in(0,1)$, denote by $V_+(w, \varrho)\equiv
\{y\in\Re^d|(y,w)_{\Re^d}\geq \varrho\|y\|_{\Re^d}\}$ the
one-sided cone with the axis $\langle w\rangle\equiv\{tw,
t\in\Re\}$, and by $V(w, \varrho)\equiv
\{y\in\Re^d||(y,w)_{\Re^d}|\geq \varrho\|y\|_{\Re^d}\}$ the
two-sided cone with the same axis.

\begin{prop}\label{p43} Suppose that the following two conditions
hold true.

\noindent 1. For every $w\in S^d,$ there exists $\varrho\in(0,1),$
such that, for  every $\delta>0$, $$
 \Pi(V(w,
\varrho)\cap\{u|\|u\|\leq \delta\})>0. $$

\noindent 2. For some point $\bx$, there exists its neighborhood
$O_\bx$ such that

2a. $c(x,u)=\chi(x)u+\delta(x,u)$, $x\in O_\bx$, and
$$
\|\delta(\bx,u)\|+\|\nabla_x\delta(\bx,u)\|=o(\|u\|),\quad
\|u\|\to 0;$$

2b. the functions $\chi(\cdot)$ and $\tilde a(\cdot)$  belong to
$C^1(O_\bx, \Re^{m\times d})$ and $C^1(O_\bx, \Re^{m})$
correspondingly, and  satisfy the following joint non-degeneracy
condition:
$$
\mathrm{rank}\,\Bigl[ \nabla \tilde a(\bx)\chi(\bx)-\nabla
\chi(\bx)\tilde a(\bx)\Bigr]=m.
$$

Then condition \tN holds true with this $x_*$ and arbitrary
$t_*>0$.
\end{prop}

\begin{rem} Condition 2b. is formulated in the case \textbf{A}. In
the case \textbf{B}, it should be replaced by the condition $\det
\nabla a(\bx)\not=0$, and, in this case, condition 2a. trivially
holds true with $\chi(x)\equiv I_{\Re^m}$.
\end{rem}

\demo We use Proposition \ref{p42}. Denote $\Bigl[ \nabla \tilde
a(\bx)\chi(\bx)-\nabla \chi(\bx)\tilde a(\bx)\Bigr]=A.$ It follows
from the condition 2 and explicit formula for $\hat \Delta$, that
\be\label{47} \hat\Delta(x_*, u)=Au+o(\|u\|),\quad \|u\|\to 0. \ee
Let $v\in S^m, \eps>0$ be fixed. Consider the linear subspace
$L_v=\{u\in \Re^d| Au\perp v\}=A^*<v>^\perp$ ($A^*$ is the adjoint
matrix for $A$). This subspace is proper, due to condition
$\mathrm{rank}\, A=m$. Take $w\in S^d$ such that $w\perp L_v,$
then, for any $\varrho\in (0,1)$, $V(w,\varrho)\cap
L_v=\emptyset$, and, therefore, there exists $c=c(v,\varrho)>0$
such that
$$
 |(Au, v)_{\Re^m}|\geq c\|u\|,
\quad u\in V(w,\varrho).
$$
 This, together with (\ref{47}),
provides that \be\label{48} |(\hat \Delta(x_*, u), v)_{\Re^m}|\geq
c\|u\|+o(\|u\|),\quad u\in V(w,\varrho), \|u\|\to 0. \ee
 Take $\varrho$ from the condition 1 of the
Proposition, and $\delta_*=\eps\cdot [\psi_*(x_*)]^{-1}$ ($\psi_*$
is given in the condition (\ref{10})). Then, for every $\delta\in
(0,\delta_*)$, the measure $\Pi$ of the  set $V(w,\varrho)\cap
\{\|u\|\leq \delta\}$ is positive, and, on this set, $\|c(x_*,
u)\|<\eps$.  On the other hand  (\ref{48}) implies that, for
$\delta$ small enough,
$$
(\hat \Delta(x_*, u), v)_{\Re^m}\not=0\hbox{ on the set }
V(w,\varrho)\cap \{\|u\|\leq \delta\}.
$$
The proposition is proved.

\subsection{Condition S}

One possible way to provide that condition \textbf{S} holds true
is to use general support theorems for the distribution of the  of
solution to SDE with a jump noise. For instance, Theorem I
\cite{t_simon_support} provides, in the case \textbf{A}, the
following result.

\begin{prop} Consider $U$, the set of sequences $\{(t_n, u_n), n\geq 1\}$,
where $\{t_n\}\subset \ax$ is a strictly increasing sequence with
$\lim t_n=+\infty$, and $\{u_n\}\subset \mathrm{supp}\,\Pi$ is
arbitrary. Suppose that, for any given $R,T\in \ax$, for every $x$
with $\|x\|\leq R$ and $\eps>0$ there exists a sequence $\{(t_n,
u_n)\}\in U$ such that the solution to the equation
$$
Z(t) = x + \int_0^t\tilde a(Z(s))\,ds + \sum_{t_n\leq t}
c(Z(t_n-), u_n),\quad t\in \ax,
$$
satisfies the condition $\|Z(T)-x_*\|<\eps$.

Then condition \textbf{S} holds true.
\end{prop}

Another possibility is to give some straightforward conditions,
that seem to be more suitable in a certain concrete cases. Let us
formulate, without a detailed proof, one condition of  such a
type. Note that, unlike the previous Proposition, the next one
does not require moment restriction on the L\'evy measure of the
noise, and is formulated for the both cases \textbf{A} and
\textbf{B} simultaneously. Denote, for any $x\in \Re^m$,
$\Pi_x(\cdot)=\Pi(u\in \Re^d| c(x,u)\in \cdot)$.

\begin{prop}
Suppose that, for every $x\in \Re^m, v\in S^m$, there exists
$\varrho\in (0,1)$ such that, for any $\delta>0$, \be\label{49}
\Pi_x(V_+(v, \varrho)\cap\{\|y\|\leq \delta\})>0. \ee

Then $y\in \mathrm{supp}\, P^t_x$ for every $x,y\in \Re^m,$ $t>0$,
and therefore condition \textbf{S} holds true.
\end{prop}

{\it Sketch of the proof}.  Take $v={y-x\over \|y-x\|}$,
$\varrho\in (0,1)$ from the condition (\ref{49}) for the given $x$
and $v$, and $\delta_*={1\over 2} \|y-x\|$. Then there exist
$\delta_1,\delta_2,\gamma>0$ such that
$0<\delta_1<\delta_2<\delta_*$,
$$
\Pi_x(V_+(v, \varrho)\cap\{\|y\|\in[\delta_1,\delta_2]\})>0,
$$
and
$$
\|(x+c)-y\|\leq \|x-y\|-\gamma, \quad c\in V_+(v,
\varrho)\cap\{\|y\|\in[\delta_1,\delta_2]\}.
$$
Then arguments, analogous to those given in the proof of
Proposition \ref{p42} allows one to conclude that, for any two
points $x\not=y$, there exist $\gamma>0$ and  $t>0$ such that, for
any $s\in (0,t)$, \be\label{410} \Ps_x(\|X(s)-y\|<
\|x-y\|-{\gamma\over 2})>0. \ee Let $\eps\in (0, \|x-y\|)$ be
given, then, since the process $X$ is Feller, one can conclude
from (\ref{410}) that there exist $t_\eps>0$, $\gamma_\eps>0$ such
that, for any $t\leq t_\eps$,
$$
p_t\equiv\inf_{z: \|z-y\|\in[\eps, \|x-y\|]}\Ps_z(\|X(s)-y\|<
\|z-y\|-\gamma_\eps)>0.
$$
This implies that, for any $x\not=y $ and $\eps>0$, for any $t\leq
t_\eps$,
$$
\Ps_x(X(t)\in B(y,\eps))\geq [p_{t\over N}]^N>0, \quad \hbox{where
}N=\Big[{\|x-y\|\over \gamma_\eps}\Big]+1.
$$
Via the Markov property of the process $X$, this implies the
statement of the Proposition.

\subsection{One-dimensional case. Proof of Proposition \ref{p01}.} In the case $m=1$,
the sufficient conditions given in the previous subsections can be
made more precise. For instance, the following version of
Proposition \ref{p42} holds true.

\begin{prop}\label{p45} Let $m=1$ and suppose that there exists $x_*\in \Re$ such that
\be\label{m1} \Pi\Big(u\in \Theta_{x_*}\Big| \hat \Delta(x_*,u)
\not=0 \Big)>0. \ee Then condition \textbf{N} holds true for any
$t_*>0$.
\end{prop}

We omit the proof, since it is  totally analogous to the one of
Proposition \ref{p42}, except one point, that causes the
difference between conditions (\ref{m}) and (\ref{m1}). For $m=1$,
we have to apply estimate (\ref{431}) only once, for the jump
moment $\tau_1^n$. This means that we do not have to control the
position of the process after the jump at this moment, and thus,
when $m=1$, the limitation involving $\eps$ can be removed from
(\ref{m}).

Now, let us prove Proposition \ref{p01}, formulated in the
Introduction.  Condition \textbf{R} is provided by Proposition
\ref{p41}. Let us proceed with the conditions \textbf{N} and
\textbf{S}. Since $\Pi(\Re\backslash \{0\})>0$, either
$\Pi((-\infty,0))$ or $\Pi((0,\infty))$ is non-zero. Let, for
instance, $\Pi((0,\infty))>0$. Take $R$ large enough for
$\sup_{x>R}{a(x)\over x}<0$, then $y\in \mathrm{supp}\, P^t_x$ for
any $y>R$ and any $x\in \Re, t>0$. This follows from Theorem I
\cite{t_simon_support} in the case $\int_{\Re}|u|\Pi(du)<+\infty$,
and from Theorem 3 \cite{t_simon_support_onedim} in the case
$\int_{\Re}|u|\Pi(du)=+\infty$. This provides that \textbf{S}
holds true with arbitrary $t>0$ and $x_*>R$. In order to provide
\textbf{N} for some $x_*>R$,  let us use  Proposition \ref{p45}.
In the case of additive noise, $\hat \Delta(x,u)=a(x+u)-a(x)$.
Therefore, if $\Pi(\Re\backslash \{0\})>0$ and (\ref{m1}) fails
for every $x_*>R$, then there exists a sequence $\{x_n\}$ with
$|x_n|\to +\infty$ such that $a(x_{n+1})=a(x_n)$. This, however,
contradicts the condition $\lim\sup\limits_{|x|\to+\infty}
{a(x)\over x}<0$. The proposition is proved.

\section{Counterexamples}

We have seen that three basic conditions
\textbf{R},\textbf{N},\textbf{S} imply exponential estimates
(\ref{111}), (\ref{12}). In this section we give counterexamples
that show that, as soon as any of these conditions is removed, the
solution to (\ref{01}) may fail to be ergodic (i.e., to possess a
unique invariant distribution $\mu_{inv} \in \Pf$).

The cases, when conditions \textbf{R} or \textbf{S} are missed,
are quite standard and simple, thus we just outline the
corresponding examples.

\begin{ex} Let $m=d=1$, $c(x,u)=u, \Pi=2\delta_1+\delta_{-1}$ and $a(x)\in
C^1(\Re)$ is such that $a(x)=-c, |x|\geq 1, $ with $c\in(0,1)$.
Both conditions \textbf{N} and \textbf{S} hold true here (this can
be  provided by the arguments from subsection 4.4), but \textbf{R}
fails. The law of large numbers provides that, for every $x>1$,
$$
\Ps_x(\lim_{t\to+\infty}X(t)=+\infty, \inf_{t\in \ax} X(t)> 1)>0.
$$
This implies that $X$ does not have any invariant probability
measure.
\end{ex}

\begin{ex} Let $m=d=1$, $a(x)\in C^1(\Re)$ be such that $a(x)=-x,
|x|\geq 2$ and $a(x)=0, |x|\leq 1$. Let also $\Pi=\delta_1$ and
$c(\cdot,1)\in C^1(\Re)$ be bounded and such that $c(x,
1)=\mathrm{sign}\, x , |x|\geq 2$ and $x\cdot c(x,1)\geq 0, x\in
\Re$. Then condition \textbf{R} holds true, and  \textbf{N} holds
true for any $x_*$ with $|x_*|>2$ (Proposition \ref{p45}).
Condition \textbf{S} fails: starting from any set
$A_+=[1,+\infty)$ or $A_-(-\infty, -1]$, the process $X$ remains
in this set with the probability 1. Therefore, there exist at
least two different invariant measures for $X$, supported by these
sets.
\end{ex}

The last example is more non-trivial, and is concerned with the
case where \textbf{R},\textbf{S} hold true while \textbf{N} does
not.

\begin{ex} Let us start with an auxiliary construction. Consider
the unit circle $C\equiv {1\over 2\pi}S^2$ on the plane $\Re^2$,
and define  the discrete time Markov process $Z$ on $C$ by its
transition probability
$$
Q(z,\cdot)=(1-3p)\delta_{3z}(\cdot)+p\left[\delta_{z\over
3}(\cdot)+\delta_{z+1\over 3}(\cdot)+\delta_{z+2\over
3}(\cdot)\right], \quad z\in C,
$$
where $p\in(0,{1\over 6})$ is given, and every arithmetic
operation on $C$ is defined as the same operation on $[0,1)\cong
C$ modulo 1. If $Z_0=z$ is any point from $[0,1)\cong C$, then
there exists a non-zero probabilities for $Z_n$ to  be equal to
each point of the type
$$
3^{-n}y+\sum_{k=1}^{n}a_k3^{-k}, \quad a_k\in\{0,1,2\}, k=1,\dots,
n,
$$
and therefore, the set $\bigcup_{n\in \NN}\mathrm{supp}\, P(Z_n\in
\cdot|Z_0=z)$ is dense in $C$, i.e. the process $Z$ is
topologically irreducible. Now let us show that $Z$ possesses at
least two different invariant measures (in fact, the set of
invariant measures here is much larger).

Consider, together with $Z$, the sequence $T_n$ defined by
$$
T_0=0, \quad T_{n+1}=\begin{cases}(T_n-1)\vee 0,&
Z_{n+1}=3Z_{n},\\
T_n+1, &\hbox{otherwise},
\end{cases}\quad n\geq 0.
$$
Then $\{T_n\}$ is a birth-and-death Markov chain with
probabilities of birth equal to $b_k\equiv b=3p$ and probabilities
of death equal to $d_k\equiv d=1-3p$. We have that $d>b$ since
$p<{1\over 6}$,  and therefore this chain is ergodic, that means
that for any given $\eps>0$ there exists $L_\eps\in \NN$ such that
\be\label{51} \sup_{n\geq 0}P(T_n\geq L_\eps)<\eps.\ee Take
$Z_0=0$ and consider some weak limit point $\mu_*^0$ for the
sequence of Khasminskii's averages $$ {1\over N}\sum_{n\leq
N}P(Z_n\in \cdot |Z_0=0)
$$ (see the proof of Corollary \ref{c34}). By the construction,
  any digit in the 3-adic representation for
$Z_n$,  with the number of the digit greater then $T_n$,  is equal
to $0$. This means that $ \mu_*^0(A_\eps^0)\geq 1-\eps, $ where
$$
A_{\eps}^0=\{z\in[0,1)|\hbox{ all 3-adic digits  for $z$, with the
number of the digit $\geq L_\eps$,  are equal to $0$}\}
$$
(the inequality holds true since the set $A_{\eps}^0$ is closed).
Therefore, $\mu_*^0(A^0)=1,$ where
$$
A^0=\{z\in[0,1)|\hbox{ all 3-adic digits  for $z$, except some
finite number of the digits, are equal to $0$}\}.
$$
Analogously, if $Z_0={1\over 2},$ and $\mu_*^{1\over2}$ is any
weak limit point for the sequence of Khasminskii's averages
$${1\over N}\sum_{n\leq N}P(Z_n\in \cdot |Z_0={1\over 2}),$$
 then $\mu_*^{1\over2}(A^1)=1,$ where
$$
A^1=\{x\in[0,1)|\hbox{ all 3-adic digits  for $x$, except some
finite number of the digits, are equal to $1$}\},
$$
and $A^0\cap A^1=\emptyset$. This means that $\mu_*^0$ and
$\mu_*^{1\over2}$ are mutually singular invariant measures for
$Z$.

Now, let us proceed with the construction of the process. Put
$m=2, d=2$ and $c(x,u)=c_1(x,u_1)+c_2(x,u_2), x\in \Re^2,
u=(u_1,u_2)\in \Re^2.$ Let $\Pi=\Pi_1\times\Pi_2$ with
$\Pi_1=\delta_1$,
$\Pi_2=(1-3p)\delta_1+p(\delta_2+\delta_3+\delta_4)$. Let the part
$c_1$ to give the radial part of the jump noise:
$$
c_1(x,1)={b(x)\over \|x\|}\cdot x, \quad x\in \Re^2,
$$
where $b\in C^1(\Re^2,\Re)$ is such that $b(x)=0$ for  $\|x\|\leq
1$, $b(x)>0$ for $\|x\|\in (1,2)$ and $b(x)=1$ for $\|x\|\geq 2$.
For $\|x\|\geq 1$, let the part $c_2$ to define the "rotational"
part of the noise: if $x$ is written in the polar coordinates as
$(r,\theta)$, then $x+c_2(x,i), i=1,\dots, 4$ has the following
polar representation:
$$
\begin{cases}(r,3\theta),&  i=1\\
(r,{\theta+2\pi(i-2)\over 3}),& i=2,3,4
\end{cases}.
$$
For $\|x\|\geq 1$, let the functions $c_2(\cdot ,i),i=1,\dots,4$
be defined in an arbitrary way, such that $c_2(\cdot ,i)\in
C^1(\Re^2, \Re^2), \|x+c(x,i)\|\geq 1, x\in\Re^2, i=1,\dots,4.$
The drift coefficient let be equal to $a(x)=-b(x)\cdot x, x\in
\Re^m$.

By the construction, condition \textbf{R} holds true (Proposition
\ref{p41}) and condition \textbf{S} holds true for every $x_*$
with $\|x^*\|\geq 1$ and every $t>0$ (Theorem I
\cite{t_simon_support}). Let us show that, however, there exist
two different invariant measures for $X$. If $X(0)=x$ is such that
$\|x\|\geq 1$, then the processes $R(\cdot)=\|X(\cdot)\|$ and
$\Theta(\cdot)={X(\cdot)\over \|X(\cdot)\|}$ are independent
(w.r.t. $\Ps_x$) Markov processes. The first process possesses at
least one invariant measure $\kap, $ supported by $[1,+\infty)$
(see Corollary \ref{c34}).
 The second one is the pure
jump Markov process with the total intensity of the jump equal, at
every point, to $(1-3p)+p+p+p=1$. Its embedded Markov chain
coincides, up to the scaling parameter $2\pi$, with the chain $Z$
considered before. Therefore, this process possesses at least two
different invariant measures $\chi_1,\chi_2$ on $S^2$. Thus, the
process $X(\cdot)$ possesses at least two different invariant
measures $\mu_1=\kap\times\chi_1,\mu_2=\kap\times\chi_2,$
supported by $[1,+\infty)\times S^2=\{x|\|x\|\geq 1\}$.

This example shows that  the topological irreducibility condition
\textbf{S}, together with the recurrence condition \textbf{R},  is
not strong enough to produce ergodicity of the solution to SDE
with a jump noise. In order to produce ergodicity, some kind of
"smoothing" condition, like non-degeneracy condition \textbf{N} in
our settings, is needed additionally.
\end{ex}

\end{document}